\documentclass[11pt]{amsart}
\usepackage{amsthm}

\usepackage{mathrsfs}
\usepackage{mathtools}
\usepackage{slashed}
\usepackage{tikz-cd}

\usepackage{amssymb}
\usepackage{amsfonts}
\usepackage{amsmath}
\usepackage{graphicx}
\usepackage{xcolor}
\usepackage{amsmath}
\usepackage{mathrsfs}
\usepackage{stmaryrd}
\usepackage{epsfig,color}
\usepackage{blindtext}
\usepackage{enumerate}
\usepackage[colorlinks=true,linkcolor=red,citecolor=blue]{hyperref}
\usepackage[noabbrev,capitalize]{cleveref}
\usepackage{geometry}
\usepackage{url}
\usepackage{nicefrac,mathtools}
\DeclareGraphicsExtensions{.pdf,.jpeg,.png}
\usepackage{epstopdf}
\usepackage{cancel} 
\usepackage[normalem]{ulem} 
\usepackage{verbatim} 
\usepackage{enumitem} 
\usepackage{color}
\usepackage[msc-links, lite, alphabetic]{amsrefs}
\usepackage{geometry}
\DeclareMathOperator{\tr}{\rm tr}
\DeclareMathOperator{\Vol}{Vol}
\DeclareMathOperator{\vol}{Vol}

\DeclareMathOperator{\Ric}{Ric}

\DeclareMathOperator{\Rm}{\rm Rm}
\DeclareMathOperator{\diam}{\mathrm{diam}}

\renewcommand{\subset}{\subseteq}

\newcommand{\eps}{\varepsilon}

\newcommand{\N}{\mathbb{N}}

\newcommand{\R}{\mathbb{R}}
\renewcommand{\subset}{\subseteq}
\newcommand{\defeq}{\mathrel{\mathop:}=}

\newcommand{\haus}{\mathcal{H}}

\newcommand{\dist}{\mathsf{d}}

\def\XXint#1#2#3{{\setbox0=\hbox{$#1{#2#3}{\int}$ }
\vcenter{\hbox{$#2#3$ }}\kern-.6\wd0}}

\def\sideremark#1{\ifvmode\leavevmode\fi\vadjust{\vbox to0pt{\vss
 \hbox to 0pt{\hskip\hsize\hskip1em
 \vbox{\hsize3cm\tiny\raggedright\pretolerance10000
 \noindent #1\hfill}\hss}\vbox to8pt{\vfil}\vss}}}

\newtheorem{theorem}{Theorem}[section]
\newtheorem{mainthm}{Theorem}

\newtheorem{proposition}[theorem]{Proposition}
\newtheorem{lemma}[theorem]{Lemma}
\newtheorem{corollary}[theorem]{Corollary}

\theoremstyle{definition}
\newtheorem{definition}[theorem]{Definition}

\newtheorem{conjecture}[theorem]{Conjecture}
\newtheorem{remark}[theorem]{Remark}

\numberwithin{equation}{section}

\allowdisplaybreaks[4]

\usepackage{esint}
\allowdisplaybreaks[2]
\makeatletter
\@namedef{subjclassname@2020}{\textup{2020} Mathematics Subject Classification}
\makeatother
\hypersetup{
 colorlinks=true,
 allcolors=black,
 pdftitle={Uniqueness of the asymptotic limits for Ricci-flat manifolds with linear volume growth},
 pdfauthor={Zetian Yan and Xingyu Zhu},
 pdfsubject={ Ricci flat, asymptotic limits, uniqueness, exponential convergence}}
\let\epsilon\eps
\newcommand{\sliceball}[1]{B_{\mathscr S^{m,\alpha}}(g_0,#1)}
\newcommand{\kernelball}[1]{B_{K_0}(0,#1)}

\begin{document}

\title[Uniqueness of cross sections]{Uniqueness of the asymptotic limits
for Ricci-flat manifolds with linear volume growth}

\author[Z. Yan]{Zetian Yan}
\address[Z. Yan]{Department of Mathematics \\ The Chinese University of Hong Kong \\ Shatin, N.T. \\ Hong Kong}
\email{zetianyan@cuhk.edu.hk}

\author[X. Zhu]{Xingyu Zhu}
\address[X. Zhu]{Department of Mathematics \\ Michigan State University \\
East Lansing \\ MI 48824 \\ USA}
\email{zhuxing3@msu.edu}


\begin{abstract}
Under natural assumptions on curvature and cross section, we establish the uniqueness of asymptotic limits and the exponential convergence rate for complete noncollapsed Ricci flat manifolds with linear volume growth, which are known to only admit cylindrical asymptotic limits. In particular, our results show that all asymptotically cylindrical Calabi--Yau manifolds converge exponentially to their asymptotic limits, thereby answering affirmatively a question by Haskins--Hein--Nordstr\"om. In dimension four, these assumptions hold automatically, yielding unconditional uniqueness and convergence thus our result strengthens those of Chen--Chen, who proved exponential convergence to its asymptotic limit space for ALH instantons. Compared to cones, the greatest challenge is to bound the multiplicity two zero modes and oscillatory negative eigenvalue modes, which are new phenomena for cylinders. To overcome the new challenge, we develop a new technical tool, the constant mean curvature foliation. 

\end{abstract}

\maketitle

\section{Introduction}

\subsection{Statement of the theorem}
 In this paper we study the uniqueness of asymptotic limits and the rate of convergence to them for Ricci-flat manifolds with linear volume growth. The motivation comes from a comparison with tangent cones at infinity for Euclidean volume growth. Let $(M,g)$ be a complete noncompact Riemannian $n$-manifold with nonnegative Ricci curvature. Gromov's precompactness theorem implies that any divergent \emph{translation sequence} $(M,g,p_i)$, meaning that $\dist_g(p_i,x)\to\infty$ for every fixed $x\in M$, has a subsequence that converges in the pointed Gromov--Hausdorff (pGH) sense to a Ricci limit space, which we call an \emph{asymptotic limit}. Calabi and Yau independently showed that the minimal order of volume growth is {linear}, see \cite{yaulinearvolumegrowth}. We therefore say that $M$ has linear (or minimal) volume growth if, for some, hence every, base point $x\in M$,
\begin{equation}\label{eq:linear-volume-growth}
 \limsup_{r\to \infty}\frac{\vol_g(B_r(x))}{r}<\infty.        
\end{equation} 
In addition, we say $M$ is \emph{noncollapsed} if
\begin{equation}\label{eq:noncollapse}
    \text{there exists $v>0$ such that $\vol_g(B_1(x))>v$ for all $x\in M$.}
\end{equation}
 Building on earlier work of Sormani \cites{SormaniMiniVol,SormaniSublinear}, the second author proved that under these assumptions every asymptotic limit is
a product $\mathbb R\times N$ with compact connected cross section \cite{Zhu2025}*{Theorem 1.2}. The cross sections may be singular and may depend on the sequence.
 It is a basic question whether asymptotic limits are unique. In general, they need not be. The cross sections may be non-isometric \cite{SormaniMiniVol}*{Example 27} or even non-homeomorphic \cite{Zhu2025}*{Theorem 5.1}, by constructions analogous to those for tangent cones at infinity \cite{CN11}. For manifolds asymptotic to a cylinder in a stronger sense, the rate of convergence is also of great interest, for instance in the study of ALH gravitational instantons \cite{ChenChen}, Calabi--Yau manifolds \cites{HaskinsHeinNordstrom,SunZhang}, and $G_2$-manifolds \cite{NordstromG2}.

In this paper we study the uniqueness and convergence rate to asymptotic limits inspired by the study of asymptotic cones by Cheeger--Tian \cite{Cheeger-Tian1994} and Colding--Minicozzi \cite{ColdingMinicozziUniqueness}. We introduce a new constant-mean-curvature foliation technique to overcome new challenges in the cylindrical setting.

At a Ricci flat metric $g$ the linearization of the Ricci tensor is
\begin{equation}
 \label{eq:deform}
 2D\Ric_g(h)=\Delta_g^Lh-2\delta_g^*\delta_gh-\nabla^2\tr_gh,
 \qquad \Delta_g^L=\nabla^*\nabla-2\Rm_g.
\end{equation}
Here,
\[
    \Delta_g^Lh=\nabla_g^*\nabla_gh-2\Rm_g(h),
    \qquad
    \delta_g^*\omega=\frac12\mathcal L_{\omega^\sharp}g.
\]
 The operator $\Delta_g^L$ is referred to as the \emph{Lichnerowicz Laplacian}. We introduce the following crucial regularity assumption.
\begin{definition}[\cite{haslhoferPhD}*{Definition 2.13}, \cite{Besse}*{Section 12}]\label{def:int}
A closed Ricci flat manifold $(N,g)$ is \emph{integrable} if every smooth symmetric tensor satisfying
\begin{equation}
 \label{eq:linearized}
 D\Ric_g(h)=0
\end{equation}
is the velocity at $t=0$ of a smooth family of Ricci flat metrics $g_t$, $|t|<\varepsilon$, with $g_0=g$. 
\end{definition}
 Now we are in a position to state our main theorem. It suffices to consider the case of one end. If $M$ has two ends then it splits as $\R\times N$. The linear volume growth assumption~\eqref{eq:linear-volume-growth} forces \(N\) to be compact. Then the uniqueness and convergence speed estimates trivially hold. We therefore assume throughout this paper that \(M\) has only one end.
\begin{mainthm}\label{thm:main}
Let $(M^n,g)$ be a complete noncompact Ricci flat $n$-manifold that is noncollapsed as in \eqref{eq:noncollapse} and has linear volume growth as in \eqref{eq:linear-volume-growth} and has one end. Assume in addition:
\begin{enumerate}[label=R.\arabic*]
 \item\label{item:integral} There is a constant $\Lambda$ such that
 for every $p\in M$,
 \begin{equation*}
  \int_{B_1(p)} |\mathrm{Rm}|^{n/2}\,d\operatorname{vol}_g
  \le \Lambda.
 \end{equation*}
 \item\label{item:limit} There is an asymptotic limit
 $(\bar N\defeq\R\times N,\bar g=dr^2+g_N)$ whose cross section $N$
 is integrable.
\end{enumerate}
Then the asymptotic cylinder is unique. There exist a compact set $K$, a number $R$, and a diffeomorphism $\Phi:(R,\infty)\times N\longrightarrow M\setminus K$ such that, for every $k\geq0$ and $0<\alpha<1$,
\begin{equation}
 \label{eq:speed}
 \|\Phi^*g-(dr^2+g_N)\|_{C^{k,\alpha}([r,r+1]\times N)}
 =O(e^{-\beta r})
\end{equation}
for some $\beta>0$. In fact the optimal $\beta$ is the square root of the first positive eigenvalue of $\Delta^L_g$.
\end{mainthm}

\begin{remark}
In the previous version of this paper, the main theorem has an extra assumption that $\Delta^L_g$ is nonnegative. This assumption is now removed. However, to the authors best knowledge, there is no example of closed Ricci flat manifold whose Lichnerowicz Laplacian has negative eigenvalues. Dai--Wei--Wang \cite{DaiWangWei} showed that any closed Riemannian manifold whose universal cover admits a nonzero parallel spinor, which includes all known examples, has nonnegative Lichnerowicz Laplacian.
\end{remark}

We will show in Proposition~\ref{prop:local-gauge} that \eqref{item:integral} implies that every asymptotic limit of $M$ is smooth. Consequently, every limit cross section $N$ is a closed Ricci flat manifold and \eqref{item:limit} is well-defined. If smoothness of all asymptotic limits is known a priori, then assumption \eqref{item:integral} is unnecessary. In particular, Theorem~\ref{thm:main} answers affirmatively a question of Haskins--Hein--Nordstr\"om \cite{HaskinsHeinNordstrom}*{concluding remarks}, namely whether every asymptotically cylindrical Calabi--Yau manifold converges exponentially to its asymptotic limit. In their setting, the cross sections are all integrable. See \cite{HaskinsHeinNordstrom}*{concluding remarks}. In contrast, for asymptotically conical Calabi--Yau manifolds, the cross sections may not always be integrable, so \cite{Cheeger-Tian1994} does not by itself establish polynomial convergence in this generality. Under Euclidean volume growth and quadratic curvature decay, polynomial convergence is established by \cite{SunZhang}*{Theorem 1.2}.

Fortunately, in the very special case $n=4$, thanks to the codimension $4$ theorem of Cheeger--Naber \cite{CN15_codim4}, all extra assumptions are not needed. We get a complete description of the uniqueness of asymptotic limits in dimension $4$.
\begin{corollary}
      Let $(M^4,g)$ be a complete non-compact Ricci flat $4$-manifold that is noncollapsed \eqref{eq:noncollapse} and has linear volume growth \eqref{eq:linear-volume-growth}. Then the asymptotic limit of $M$ is unique and the convergence on each end, defined as in \eqref{eq:speed}, is exponential.
\end{corollary}
 Even though we no longer need the local $L^{n/2}$-curvature bound \eqref{item:integral} when $n=4$, it in fact holds as shown by Jiang--Naber \cite{JiangNaber}. This corollary strengthens the result of Chen--Chen \cite{ChenChen}*{Theorem 4.18}, who proved exponential convergence for ALH gravitational instantons under an a priori polynomial convergence assumption. Our result applies to a larger class of manifolds and requires no a priori convergence rate.

\subsection{Background and Motivation}\label{sec:motivation}
We explain why asymptotic limits for linear (minimal) growth are natural analogues of tangent cones at infinity for Euclidean (maximal) volume growth. Recall that $(M,g)$ has nonnegative Ricci curvature. We say $M$ has Euclidean volume growth if the monotone decreasing volume ratio has positive limit at infinity for some (hence any) base point $x\in M$:
\begin{equation}\label{eq:max-volume-growth}
   \lim_{r\to\infty}\frac{\vol_g(B_r(x))}{r^n}>0.
\end{equation} 
Note that Euclidean volume growth automatically implies noncollapsedness. 

Given any scaling sequence $(M,r_i^{-2}g)$ with $r_i\to \infty$, a subsequence pointed GH-converges to a Ricci limit space, called a \emph{tangent cone at infinity} (or an asymptotic cone). When $M$ has Euclidean volume growth, every asymptotic cone is a metric cone $(C(N)=[0,\infty)\times N, dr^2+r^2g_N)$ \cite{Cheeger-Colding97I}, with the cross section $N$ depending on the scaling sequence. This theorem is known as the volume cone to metric cone theorem. It can be viewed as the (almost) rigidity theorem when the monotone quantity ${\vol_g(B_r(x))}/{r^n}$, equivalently the surface area quotient ${\haus^{n-1}(\partial B_r(x))}/{r^{n-1}}$, (almost) becomes a constant.

For asymptotic limits for linear growth, there is a parallel picture. The cylinder structure of asymptotic limits, just like the cone structure of tangent cones at infinity for Euclidean growth, comes from the rigidity when some monotone quantity becomes a constant. This monotone quantity is defined by a \emph{Busemann function}. Let $\gamma:[0,\infty)\to M$ be a ray, the Busemann function $b_\gamma$ associated to $\gamma$ is defined to be 
\[
b_\gamma(x)=\lim_{t\to\infty} t-\dist_g(x,\gamma(t)).
\]
It is shown by Sormani \cite{SormaniMiniVol}*{Theorem 19, Lemma 20} that when $M$ has linear volume growth, then $b_\gamma$ is proper and $r\to \haus^{n-1}(\{b_\gamma=r\})$ is monotone increasing. It is then shown by the second author \cite{Zhu2025}*{Proposition 3.5} that with the extra assumption that $M$ is noncollapsed, $\diam(\{b_\gamma=r\})$ is uniformly bounded. Combining the uniform bound on the diameter of level sets and the monotonicity, we see that the cylinder structure of asymptotic limits indeed follows from the almost rigidity theorem when the monotone quantity $\haus^{n-1}(\{b= r\})$ almost becomes a constant by Sormani \cite{SormaniSublinear}*{Theorem 34}, after Cheeger--Colding \cite{CheegerColding96}.

We therefore believe \cite{ColdingMinicozziUniqueness} holds for cylinders as well.
\begin{conjecture}
    Let $(M^n,g)$ be a noncollapsed Ricci flat manifold with linear volume growth. If one asymptotic limit $\R\times N$ is smooth, then the asymptotic limit is unique.
\end{conjecture}

We end this subsection by the comparison of convergence speed in the cylindrical and conical cases. It explains why exponential convergence is the correct order to expect when having polynomial convergence for cones.
\begin{remark}
    A cone metric $dt^2+t^2g$ with a change of variable $t=e^r$ is conformal to the cylinder metric $dr^2+g$. It is immediate to see that polynomial (resp. logarithmic) decay in the conical coordinates corresponds to exponential (resp. polynomial) decay in the cylindrical coordinates.   
\end{remark}

\subsection{Sketch of proof and remarks}
We give an outline of the proof to highlight the new ideas involved in the proof and the relations to the work for cones  by \cite{Cheeger-Tian1994} and \cite{ColdingMinicozziUniqueness}. 

The first new ingredient in the proof is the foliation of the end of $M$ by hypersurfaces of constant mean curvature (CMC) in Section \ref{sec:construction-gauge}. The hypersurfaces will be parametrized by $r$. Heuristically, a CMC foliation is possible as we have reference diffeomorphisms to the limit cylinders $\R\times N$ where each slice $\{r\}\times N$ is a leaf with constant mean curvature $0$. It suggests a small perturbation can produce surfaces of constant mean curvature on the end of $M$. A delicate implicit function argument is needed to rigorously carry out the construction. 

The motivation of CMC foliation is to construct a geometric monotone quantity, which is the mean curvature along leaves of the foliation. Indeed, it will be shown that if $\gamma(r)$ is the metric on the leaf with $g=u^2\, dr^2+\gamma$ and proper normalization of $u$, $H(r)$ its mean curvature, and $A(r)$ its second fundamental form, then
\begin{equation*}
 (\log\vol_\gamma N)'=H,\qquad
 H'=-\fint_N u|A|^2\,d\vol_\gamma\leq0.
\end{equation*}
This idea comes from the Green distance monotone quantity by Colding--Minicozzi. In the linear growth case, the asymptotic behavior of linear growth harmonic functions is not well explored compared to Green's function in the Euclidean growth case, so instead of constructing a similar monotone quantity by harmonic functions, we construct a geometric monotone quantity $H$.

The second ingredient is a decomposition of the metric with a divergence condition, which is carried out for leaves converging to the integrable metric $g_N\defeq g_0$ in \eqref{item:limit}. Since this is achieved by the Ebin slice theorem and integrability, we refer to this divergence condition as the fixed slice condition, which forms the main body of Section \ref{sec:fixed-slice}. More precisely, we introduce another diffeomorphism $\Phi_\phi$ constructed by a shift vector field $X$, which will be called a slice coordinate map. With this shift we can find a family of metrics $g_z$ parametrized by a symmetric $(0,2)$-tensor $z\in K_0=\ker\Delta^L_{g_0}\cap\ker\delta_{g_0}$ so that $\delta_{g_0}(\gamma-g_z)=\delta_{g_0}(g_z-g_0)=0$, $g_z=g_0+z+\eta(z)$, and $\gamma-g_z=h$. The tensor $h$ will be called the \emph{transverse} term and $z$ the \emph{kernel} term. $h$ has trivial projection onto $K_0$ hence can be further decomposed into the negative and positive spectral components $h=h_-+h_+$. The diffeomorphism $\Phi$, serving as a chart to the limit cylinder, in which the aforementioned decomposition exists, is only constructed locally along the sequence in \eqref{item:limit}. Extending this chart to the whole end will allow some estimates that hold in the chart to pass to infinity. The construction of $\Phi$, along with a criterion for extension, will be given in Proposition \ref{prop:slice-path}.

\begin{remark}
The slice coordinate map $\Phi_\phi$ and the decomposition of metrics on the cross section are a reminiscent of the divergence-free gauge and solutions of the linearized Ricci flat equation in Cheeger--Tian \cite{Cheeger-Tian1994}*{Sections 3 and 4}. 

The slice coordinate map serves similar purpose of the divergence free gauge in \cite{Cheeger-Tian1994}*{Section 3}, there is a similar step to extend this gauge to infinity.     
    
The linearized equation $D\Ric(h)=0$ on a cylinder with divergence free condition can be reduced to $(-\partial_r^2+\Delta^L_{g_0})h=0$. When solving this equation by spectral decomposition of $\Delta^L_{g_0}$, it has solutions $B+r\widetilde B$ for the zero eigenvalue, where $B$, $\widetilde B$ are transverse-traceless tensors. The linear term $r\widetilde B$ is present as the characteristic equation has $0$ as a multiplicity two solution. A constant kernel term $B$ can be absorbed into a Ricci flat metric near the reference Ricci flat metric using integrability, whereas the term $rB$ cannot be removed by a fixed choice of a nearby Ricci flat metric. Moreover, $h_-$ represents oscillatory solutions with no decay or growth. The oscillatory terms cannot be directly controlled by the three circle estimates used in the conical argument of \cite{Cheeger-Tian1994}*{Section 5}. We will see how the CMC foliation controls them.
\end{remark}

The two ingredients are put together in Section \ref{sec:uniqueness} by the contracted Gauss equation on a leaf, which reads
\[
R_\gamma=H^2-|A|^2
\]
where $R_\gamma$ is the scalar curvature on $N$ with respect to the metric $\gamma$ on the leaf. The integral of this identity relates $H$ and $A$ to the Einstein--Hilbert functional $E(g)=\int_N R_g\,dV_g$. Its second variation for $h$ at $g_0$ has the Lichnerowicz Laplacian as the principal term. More explicitly,
\begin{equation*}
 D^2E_{g_0}(h,h)
 =-\tfrac12\langle\Delta^L_{g_0}h,h\rangle_{L^2} +\tfrac12\|d\tr_{g_0}h\|_{L^2}^2.
\end{equation*}
It is well known that Ricci flat metrics are critical points of $E$, so the estimates of $E$ by Taylor expansion will be directly related to $D^2 E$. This is how the eigenvalues of $\Delta^L_g$ come into play. Inspired by the work of Bamler and Lai \cite{BamlerLai2025}, we establish the key estimates
\begin{equation}
 \int_N R_\gamma\,dV_\gamma
 \ge c\|h_-\|_{L^2(N)}^2-C\|h_+\|_{\mathsf{H}^1(N)}^2
\end{equation}
and
\begin{equation}
 |z'|+\|A\|_{L^2(N)}+\|h_-\|_{\mathsf{H}^1(N)}
 \le C\bigl(H+\|h_+\|_{\mathsf{H}^1(N)}\bigr).
\end{equation}
The last estimate controls both the oscillation part $h_-$ and the kernel term $z$ after an integration.  

To derive the exponential decay, we need to estimate the transverse term $h_+$ and mean curvature $H$. We divide the interior of the interval $[a,T)$, on which the slice coordinates are well-defined, into subintervals $I_j$ of fixed length $L$. Let  $e_j=\|h_+\|_{\mathsf{H}^2(I_j\times N)}$,
 $H_j=\sup_{I_j}|H|$, then on each $I_j\times N$, it holds $H+\|h_+\|_{\mathsf{H}^1(N)}\le H_j+C e_j$ and we can prove for some $\theta\in (0,1)$
\begin{equation}
 e_j\leq C\bigl(\theta^je_0+\theta^{J-j}e_J\bigr) +C\varepsilon\sum_{k=0}^J\theta^{|j-k|}(H_k+e_k).
\end{equation}  
Here $J$ is the last index of the consecutive subintervals in $[a,T)$. This is the replacement of the three circle type estimates in Cheeger--Tian \cite{Cheeger-Tian1994}*{Section 5}. Once the slice coordinates can be extended to infinity, the estimates in $I_j$ can also extend to infinity, i.e. we can let $J\to\infty$. Finally, the above inequality of $e_j$ gives rise to a difference inequality of $\sup_{k\ge j}(H_k+e_k)$ and the exponential decay follows from a discrete Gronwall argument.

\begin{remark}
In the previous version of this paper, we tried to prove the main theorem by Cheeger--Tian's framework but there is a gap in the proof about absorbing the multiplicity two zero mode $rB$ into the three circle type estimates. The new method was introduced to fix the gap but also helped us to remove the assumption on the sign of Lichnerowicz Laplacian. We do not know if the original method of \cite{Cheeger-Tian1994}, inspired by \cites{SiomonAsymptotics}, still work for cylinders as Haskins--Hein--Nordst\"orm \cite{HaskinsHeinNordstrom} suggested.
\end{remark}

\subsection{Conventions and organization}
A cylinder is the product of $\mathbb R$ with a closed connected manifold. 

We write $\mathsf{H}^j(\Omega)$ for the $L^2$ Sobolev space of order $j$ on the domain $\Omega$. A prime denotes $\partial_r$.
Product norms on $I\times N$ include all mixed derivatives in the $r$ and $N$ directions. In Sections \ref{sec:fixed-slice} and \ref{sec:uniqueness}, leafwise norms and inner products use the metric $g_0=g_N$, and product norms use $dr^2+g_0$, unless a different metric is explicitly used. A norm of a leafwise field with no domain notation is taken on $N$. Constants may depend on a fixed derivative order and background neighborhood.

The sign convention for the Laplacian is $\Delta=d^*d+dd^*\ge0$. On one-forms we also write $\Delta_H$.

The paper is organized as follows.  Section~\ref{sec:construction-gauge} constructs the global CMC foliation
and derives its volume identities. Section~\ref{sec:fixed-slice} constructs the fixed slice coordinates and the decomposition into Ricci flat metrics $g_z$, kernel and transverse terms $z$ and $h$, and proves estimates about the shift vector field. Section~\ref{sec:uniqueness} uses the Gauss equation constraint to prove the transverse estimate, and proves the exponential decay after extending the slice coordinates to infinity.

\subsection*{Acknowledgements}
The authors thank Junsheng Zhang for sharing the reference \cite{HaskinsHeinNordstrom}, for many invaluable discussions, and for his careful reading of an earlier draft of this manuscript. We are also grateful to Xi-Nan Ma and Guofang Wei for their interest in this work, to Yi Lai for invaluable discussions and useful comments. Z.Y. is supported by an AMS--Simons Travel Grant and acknowledges the hospitality of the University of Science and Technology of China, where part of this research was carried out during his visit. X.Z. is supported by an AMS--Simons Travel Grant.

\subsection*{Disclosure on AI Use}
ChatGPT 5.6 Sol and 6 Astra/Codex were used substantially to refine the arguments, with the major steps and ideas, including CMC foliations, slice condition and Gauss constraint, given in the prompts by the authors. But AI handled many technical arguments. In particular, some hard estimates like Lemmas~\ref{lem:transverse-preliminary}--\ref{lem:transverse-estimate} are accomplished by AI. The ideas presented in this paper are developed by the authors partially in order to fix an essential mistake in the previous version. The authors reviewed, edited, and remain responsible for all theorem statements, proofs and citations.

\tableofcontents

\section{A global CMC foliation of the end}\label{sec:construction-gauge}
In this section, we use \eqref{eq:noncollapse}, \eqref{eq:linear-volume-growth}, and condition~\ref{item:integral} to construct a global CMC foliation of the end. Integrability in condition~\ref{item:limit} is used in Section~\ref{sec:fixed-slice}. 

Fix a unit-speed ray \(\gamma_0:[0,\infty)\to M\), and let $b(x)$ be its associated Busemann function. We will repeatedly use the following consequences of \cites{SormaniMiniVol,Zhu2025}. In particular, \cite{Zhu2025}*{Lemma 2.3 and Propositions 3.2 and 3.5}:
\begin{equation}
 \label{eq:busemann-bounds}
 \begin{gathered}
 \text{the sublevel set } \{b\leq R\}\ \text{is compact for every }R,\\
  \sup_{R\in\mathbb R}\diam b^{-1}(R)\leq D,
  \qquad
  \sup_{x\in M}d(x,\gamma_0)<\infty.
 \end{gathered}
\end{equation}
These estimates make \(b\) a radial Lipschitz coordinate on the end. Properness of the Busemann sublevel sets ensures that every bounded Busemann tube is compact, so all constructions carried out over a finite range of Busemann levels take place in a compact region. The uniform diameter bound for the level sets, together with the bounded distance from \(M\) to the reference ray \(\gamma_0\), further implies
that, for each fixed $L$, the tube $b^{-1}([R-L,R+L])$ is contained in a single pointed convergence chart for all sufficiently large $R$. We first use these facts to establish uniform smooth cylindrical structure on such tubes in the proposition below. The resulting product charts will then be used to construct local CMC families. Uniqueness on overlaps identifies these local families, and the resulting global foliation will be shown to exhaust the entire end.

\begin{proposition}\label{prop:local-gauge}
Under \eqref{eq:linear-volume-growth},  \eqref{eq:noncollapse}
and condition~\ref{item:integral}, all asymptotic limits are smooth
Ricci-flat cylinders. Their cross sections are all diffeomorphic
and have a common volume: there is a constant $V_\infty>0$ such that
\begin{equation}
 \label{eq:common-cross-section-volume}
 \vol(N_\infty,g_\infty)=V_\infty
 \quad\text{for every asymptotic cylinder }
 (\mathbb R\times N_\infty,dt^2+g_\infty).
\end{equation}
The family of metrics \(\{g_\infty\}\) on the limit cross section is precompact in the smooth topology modulo diffeomorphism. Consequently, for every \(k\geq 0\), their curvature tensors have uniformly bounded \(k\)-th covariant derivatives. Moreover, their injectivity radii are uniformly bounded below by a positive constant, and the first positive eigenvalues of their scalar Laplacians admit a uniform positive lower bound.

For every fixed \(L>2\), integer \(k\geq 0\), and
\(\alpha\in(0,1)\), there exists \(R_0=R_0(L,k,\alpha)\) such that, for
each \(R\geq R_0\), there are a closed manifold \(N_R\), a Ricci-flat
metric \(g_R\) on \(N_R\), and a diffeomorphism onto its image
\[
 \Phi_R:(-L,L)\times N_R\longrightarrow M
\]
with product metric
\[
 \bar g_R \defeq dt^2+g_R,
\]
such that
\begin{align}
 \|\Phi_R^*g-\bar g_R\|_
 {C^{k,\alpha}((-L,L)\times N_R)}
 &=o_R(1),
 \label{eq:smoothclose}\\
 \|b\circ\Phi_R-(R+t)\|_{C^0((-L,L)\times N_R)}
 &=o_R(1),
 \label{eq:levelclose}\\
 b^{-1}\bigl((R-L+1,R+L-1)\bigr)
 &\subset \operatorname{im}\Phi_R.
 \label{eq:local-cylinder-coverage}
\end{align}
We choose each \((N_R,g_R)\) from the family of cross sections of asymptotic cylinders. In particular, $\vol(N_R,g_R)=V_\infty$.
If
\[
 \widehat g_R
  \defeq \left.\Phi_R^*g\right|_{\{0\}\times N_R},
\]
then \eqref{eq:smoothclose} and the identity
\(\vol(N_R,g_R)=V_\infty\) imply
\begin{equation}
 \label{eq:approximate-cross-section-volume}
 \vol(N_R,\widehat g_R)\longrightarrow V_\infty,
 \qquad
 \frac12V_\infty
 \leq \vol(N_R,\widehat g_R)
 \leq 2V_\infty
 \quad\text{for }R\geq R_1,
\end{equation}
for some \(R_1<\infty\). All occurrences of \(o_R(1)\) in
\eqref{eq:smoothclose}--\eqref{eq:levelclose} are uniform as
\(R\to\infty\).
\end{proposition}

\begin{remark}
In the preceding proposition, neither the metric \((N_R,g_R)\) nor the chart \(\Phi_R\) is canonical, but all estimates are uniform over the resulting family. 
\end{remark}

\begin{proof}

Let \(p_i\to\infty\) be any divergent sequence. By \cite{Zhu2025}*{Theorem 1.2}, after passing to a subsequence,
\((M,g,p_i)\) converges in the pointed Gromov--Hausdorff sense to a
 limit cylinder $\mathbb R\times N_\infty$. 

Together with \eqref{item:integral}, the scale-invariant Einstein \(\varepsilon\)-regularity theorem gives smooth convergence away from a locally finite singular set. See \cites{AndersonEinsteinMetric,AndersonConvergence} and \cite{CN15_codim4}*{Section 2.2, (2.7)}. The singular set of the limit is contained in this concentration set. On the other hand, translation invariance in the \(\mathbb R\)-direction would carry any singular point of the product to an entire singular line, contradicting local finiteness. Hence the limiting cylinder  $\mathbb R\times N_\infty$ is smooth and carries a Riemannian metric $dt^2+g_\infty$. The Einstein equation in harmonic coordinates upgrades the pointed Gromov--Hausdorff convergence to smooth convergence in every finite \(C^{k,\alpha}\) topology.

By \cite{Zhu2025}*{Theorem 1.2 and Proposition 3.8(2)--(3)}, the Hausdorff measure of every limiting cross section equals the common limit of the perimeters of the Busemann sublevel sets. Since the limit cross sections are smooth, their Hausdorff measures agree with their Riemannian volumes. We define
\begin{equation*}
 V_\infty
  \defeq 
 \lim_{s\to\infty}\operatorname{Per}_g(\{b<s\})
 =
 \vol(N_\infty,g_\infty)
 >0.
\end{equation*}
This proves \eqref{eq:common-cross-section-volume}.

We next prove that the asymptotic cross sections are smoothly
precompact modulo diffeomorphism. Let
\[
 \bigl(\mathbb R\times N_\nu,dt_\nu^2+g_\nu,(0,y_\nu)\bigr),
 \qquad \nu\in\mathbb N,
\]
be an arbitrary sequence of asymptotic cylinders. Each asymptotic
cylinder inherits the uniform lower volume bounds and the local
$L^{n/2}$ curvature bounds \eqref{item:integral} from its
smooth realizations in $M$. Thus the same regularity argument applies to this sequence of asymptotic cylinders and upgrades its convergence to smooth pointed
convergence. After passing to a subsequence, we have the smooth convergence
\[
 \bigl(\mathbb R\times N_\nu,dt_\nu^2+g_\nu,(0,y_\nu)\bigr)
 \longrightarrow
 \bigl(\mathbb R\times N_\infty,dt^2+g_\infty,(0,y_\infty)\bigr).
\]

It is well known by a diagonal argument that the iterated asymptotic limit $(\R\times N_\infty,dt^2+g_\infty)$ is also an asymptotic limit of $(M,g)$. The diameter estimate in \eqref{eq:busemann-bounds}, passed to the asymptotic limits, gives
\[
 \diam(N_\nu,g_\nu)\leq D.
\]
Hence every central cross section
\(\{t_\nu=0\}\cong N_\nu\) is contained in a ball of fixed radius about
\((0,y_\nu)\). We may therefore choose the smooth convergence maps so
that their images contain the entire central cross sections. Pull back
\(t_\nu\) by these maps and denote the resulting functions by \(\widetilde t_\nu\). Since \(dt_\nu\) is a parallel unit one-form on the \(\nu\)-th product cylinder and $t_\nu(0,y_\nu)=0$, the functions \(\widetilde t_\nu\) converge smoothly on compact subsets, after passing to a further subsequence and with proper normalization, to the $\R$-coordinate function $t$ on $\R\times N_\infty$.

It follows that, for some \(\varepsilon_\nu\to0\),
\[
 \Psi_\nu^{-1}\bigl(\{t_\nu=0\}\bigr)
 \subset
 \bigl\{|t|<\varepsilon_\nu\bigr\},
\]
where \(\Psi_\nu\) denotes the corresponding convergence map. Moreover,
on a fixed neighborhood of $\{t_\nu=0\}$, $\partial_t\widetilde t_\nu\longrightarrow1$. For every fixed sufficiently small \(\varepsilon>0\), we therefore have, for all sufficiently large \(\nu\),
\[
 \widetilde t_\nu(-\varepsilon,y)<0
 <\widetilde t_\nu(\varepsilon,y),
 \qquad
 \partial_t\widetilde t_\nu>\frac12
 \quad\text{for every }y\in N_\infty.
\]
Consequently, every radial fibre meets
\(\{\widetilde t_\nu=0\}\) in exactly one point. The implicit function
theorem then gives a smooth function
\(\varphi_\nu:N_\infty\to\mathbb R\) such that
\[
 \Psi_\nu^{-1}\bigl(\{t_\nu=0\}\bigr)
 =
 \bigl\{(\varphi_\nu(y),y):y\in N_\infty\bigr\},
 \qquad
 \varphi_\nu\longrightarrow0
\]
smoothly. The resulting graph identifications
\(\iota_\nu:N_\infty\to N_\nu\) satisfy
\[
 \iota_\nu^*g_\nu\longrightarrow g_\infty
\]
in the smooth topology.

 The spectral bound follows from Ding under just a Ricci lower bound but we can provide an elementary proof for Ricci flat metrics. First, uniform bounds for all curvature derivatives and a uniform positive lower bound for the injectivity radii follow from smooth precompactness. The first positive scalar Laplace eigenvalues also admit a uniform positive lower bound. Indeed, otherwise, after passing to a subsequence, there would be eigenfunctions \(u_\nu\) satisfying
\[
 \begin{gathered}
 \int_{N_\nu}u_\nu\,dV_{g_\nu}=0,\qquad
 \int_{N_\nu}u_\nu^2\,dV_{g_\nu}=1,\\
 \int_{N_\nu}|\nabla u_\nu|_{g_\nu}^2\,dV_{g_\nu}
 =\lambda_1(N_\nu,g_\nu)\longrightarrow0.
 \end{gathered}
\]
After pulling back by \(\iota_\nu\), smooth compactness and the Rellich theorem give a strong \(L^2\) limit \(u_\infty\) on \(N_\infty\).
The vanishing Dirichlet energy makes \(u_\infty\) constant. Since \(N_\infty\) is connected and \(u_\infty\) has mean zero, it must vanish, contradicting its \(L^2\)-normalization.

We next prove that all asymptotic cross sections have the same diffeomorphism type.
In the space of pointed proper metric spaces endowed with the pointed
Gromov--Hausdorff topology, define the ray-based asymptotic limit set
\[
 \mathscr A_{\gamma_0}
  \defeq 
 \bigcap_{\tau\geq0}
 \overline{
  \bigl\{[(M,g,\gamma_0(s))]:s\geq\tau\bigr\}
 }.
\]
The map
\[
 s\longmapsto [(M,g,\gamma_0(s))]
\]
is continuous. By the pointed precompactness above, for large $\tau$ each closure defining \(\mathscr A_{\gamma_0}\) is compact. It is also connected, being the closure of the continuous image of the connected interval \([\tau,\infty)\). Since these compact connected sets are nested, their intersection \(\mathscr A_{\gamma_0}\) is nonempty and
connected.

Every element of \(\mathscr A_{\gamma_0}\) is a pointed asymptotic cylinder. Smooth compactness implies that whenever a sequence of such cylinders converges to another element of \(\mathscr A_{\gamma_0}\), their cross sections are eventually diffeomorphic to the cross section of the limit. Hence the diffeomorphism type of the cross section is locally constant on \(\mathscr A_{\gamma_0}\), and therefore constant because \(\mathscr A_{\gamma_0}\) is connected.

It remains to relate an arbitrary asymptotic limit to this ray-based
limit set. Let \(p_\nu\to\infty\) and suppose, after passing to a
subsequence, that $(M,g,p_\nu)\longrightarrow (X,d_X,p_\infty)$ in the pointed Gromov--Hausdorff topology. Properness of the Busemann sublevel sets gives $b(p_\nu)\longrightarrow\infty$.

Set
\[
 y_\nu \defeq \gamma_0\bigl(b(p_\nu)\bigr).
\]
Since \(b(\gamma_0(s))=s\), the points \(p_\nu\) and \(y_\nu\) lie in
the same Busemann level set. The diameter bound in
\eqref{eq:busemann-bounds} therefore gives
$d(p_\nu,y_\nu)\leq D$. After passing to a further subsequence, the marked points \(y_\nu\) converge in the pointed limit to some \(y_\infty\in X\), and $(M,g,y_\nu)\longrightarrow (X,d_X,y_\infty)$. Because \(y_\nu=\gamma_0(b(p_\nu))\) and \(b(p_\nu)\to\infty\), the
pointed space \((X,d_X,y_\infty)\) belongs to \(\mathscr A_{\gamma_0}\). This proves the assertion for all asymptotic limits.

We finally establish the product charts and properties \eqref{eq:smoothclose}-\eqref{eq:local-cylinder-coverage}. For \(R\) sufficiently large and
\[
 x\in b^{-1}\bigl([R-L-2,R+L+2]\bigr),
\]
the point \(\gamma_0(b(x))\) is defined and lies in the same Busemann
level set as \(x\). Hence
\[
 d\bigl(x,\gamma_0(R)\bigr)
 \leq
 d\bigl(x,\gamma_0(b(x))\bigr)
 +d\bigl(\gamma_0(b(x)),\gamma_0(R)\bigr)\leq D+|b(x)-R|
 \leq D+L+2.
\]
In particular,
\begin{equation}
 \label{eq:busemann-tube-ball}
 b^{-1}\bigl([R-L-2,R+L+2]\bigr)
 \subset
 B_{D+L+3}\bigl(\gamma_0(R)\bigr).
\end{equation}

To verify uniformity, suppose that one of the conclusions
\eqref{eq:smoothclose}--\eqref{eq:local-cylinder-coverage} fails along
a sequence \(R_\nu\to\infty\). After passing to a subsequence, choose
smooth pointed convergence charts centered at \(\gamma_0(R_\nu)\),
defined on a strictly larger product neighborhood and whose images
contain the ball on the right-hand side of
\eqref{eq:busemann-tube-ball}. Recall that \(t\) denotes the $\R$-coordinate of the limiting product cylinder. Then
\begin{equation}\label{eq:Busemann-alignment}
    b\circ\Phi_{R_\nu}-(R_\nu+t)\longrightarrow0 
\end{equation}
uniformly on the relevant compact product neighborhood.

For all sufficiently large \(\nu\), a point satisfying
\[
 |b-R_\nu|<L-1
\]
 has product coordinate
\[
 |t|
 \leq |b-R_\nu|
     +\bigl|b\circ\Phi_{R_\nu}-(R_\nu+t)\bigr|
 <L.
\]
Together with \eqref{eq:busemann-tube-ball}, this shows that
\[
 b^{-1}\bigl((R_\nu-L+1,R_\nu+L-1)\bigr)
 \subset
 \Phi_{R_\nu}\bigl((-L,L)\times N_{R_\nu}\bigr),
\]
which is \eqref{eq:local-cylinder-coverage}. Smooth pointed convergence
gives \eqref{eq:smoothclose}, while \eqref{eq:Busemann-alignment} gives \eqref{eq:levelclose}. This is a contradiction. Since the sequence \(R_\nu\to\infty\) was arbitrary, \eqref{eq:smoothclose}--\eqref{eq:local-cylinder-coverage} hold uniformly for all sufficiently large \(R\).
\end{proof}

We now use Proposition~\ref{prop:local-gauge} to construct a global CMC
foliation of the end and its normal transport coordinates.

\begin{theorem}\label{thm:cmc-foliation}
After deleting a compact subset, the end admits a smooth foliation by
closed connected constant-mean-curvature hypersurfaces. More precisely, there are a
compact domain \(K\subset M\) with smooth boundary and a
diffeomorphism
\[
 \Phi:[0,\infty)\times N
 \longrightarrow M\setminus\operatorname{int}K
\]
such that
\begin{equation}
 \label{eq:cmc-metric}
 \Phi^*g=u^2dr^2+\gamma(r),\qquad
 V(r) \defeq \vol(N,\gamma(r)),\qquad
 \frac1{V(r)}\int_Nu\,d\vol_{\gamma(r)}=1.
\end{equation}
Both \(V\) and \(u\) are bounded above and below by positive
constants.

For fixed \(L>0\), \(k\geq0\), \(\alpha\in(0,1)\), and
\(r_j\to\infty\), there are \((N_j,g_j)\) from
Proposition~\ref{prop:local-gauge} and diffeomorphisms
\(\vartheta_j:N_j\to N\). Set
\[
 \mathcal C_j(s,x) \defeq \Phi(r_j+s,\vartheta_j(x)),
 \qquad G_j \defeq \mathcal C_j^*g.
\]
Then, with \(u_j\) and \(A_j\) denoting the normal speed and second
fundamental form in these coordinates,
\begin{equation}
 \label{eq:cmc-local-product-limit}
 \|G_j-(ds^2+g_j)\|_{C^{k,\alpha}([-L,L]\times N_j)}+\|u_j-1\|_{C^{k,\alpha}([-L,L]\times N_j)}
 +\|A_j\|_{C^{k,\alpha}([-L,L]\times N_j)}
 \longrightarrow0,
\end{equation}
where all norms are computed with respect to \(ds^2+g_j\).
\end{theorem}

\begin{proof}
We break the rather long proof into several claims.

\textbf{Claim 1: There exists a family of constant-mean-curvature hypersurfaces.}

Apply Proposition~\ref{prop:local-gauge} with \(L=10\), choosing the charts diagonally in derivative order: at the $i$-th center require convergence up to order $m_i$, where $m_i\to\infty$.
This is possible by increasing the radius thresholds at each fixed order. Consequently the chosen charts converge in every fixed order.
Use centers
\[
 R_i=R_\ast+i,\qquad i\in\N,
\]
where \(R_\ast\) will be chosen sufficiently large. Write
\[
 N_i \defeq N_{R_i},\qquad
 g_i \defeq g_{R_i},\qquad
 \widehat G_i \defeq \Phi_{R_i}^*g,\qquad
 \overline G_i \defeq dt^2+g_i,
\]
and, for every fixed \(m\geq0\), set
\[
 \varepsilon_i^{(m)}
  \defeq 
 \|\widehat G_i-\overline G_i\|_
 {C^{m,\alpha}((-6,6)\times N_i)}.
\]
By Proposition~\ref{prop:local-gauge},
\(\varepsilon_i^{(m)}\to0\) as \(i\to\infty\), and
\(\sup_{i\geq0}\varepsilon_i^{(m)}\) can be made arbitrarily small by
increasing \(R_\ast\).

For \(|a|\leq5\) and
\(\varphi\in C^{2,\alpha}(N_i)\) small, let
\(\mathscr H_i^G(a,\varphi)\) be the mean curvature of the graph $t=a+\varphi(x)$, with respect to an ambient metric \(G\) near \(\overline G_i\), oriented by the normal with positive \(\partial_t\)-component and pulled back to \(N_i\). This is a perturbation of the level $\{t=a\}$ in the given chart. For the constant mean curvature equation, use the Banach spaces
\begin{equation}
 \label{eq:cmc-function-spaces}
 \mathcal X_i \defeq C^{2,\alpha}(N_i)\times\mathbb R,\qquad
 \mathcal Y_i \defeq C^{0,\alpha}(N_i)\times\mathbb R.
\end{equation}
The norms are taken with respect to $g_i$, and on $\mathcal X_i$ we write
\begin{equation}
 \label{eq:cmc-function-ball}
 B_{\mathcal X_i}(0,\rho_*) \defeq 
 \{(\varphi,c)\in\mathcal X_i:\|(\varphi,c)\|_{\mathcal X_i}<\rho_*\}.
\end{equation}
Define
\[
 \langle f\rangle_i
  \defeq 
 \frac1{\vol(N_i,g_i)}
 \int_{N_i}f\,d\vol_{g_i}
\]
and
\begin{equation}
 \label{eq:cmc-augmented-map}
 \mathscr A_{i,a}^G(\varphi,c)
  \defeq 
 \bigl(\mathscr H_i^G(a,\varphi)-c,\langle\varphi\rangle_i\bigr).
\end{equation}
Thus \(\mathscr A_{i,a}^G(\varphi,c)=0\) means that the graph has
constant mean curvature \(c\) and that its average radial displacement
vanishes.

At the product metric \(\overline G_i\), the slices
\(\{t=a\}\) are totally geodesic. With the convention
\(\Delta=\nabla^*\nabla\), the linearization of mean curvature is
\[
 D_\varphi\mathscr H_i^{\overline G_i}(a,0)w
 =\Delta_{g_i}w.
\]
Hence
\[
 L_i
  \defeq 
 D_{(\varphi,c)}
 \mathscr A_{i,a}^{\overline G_i}(0,0)
\]
is independent of \(a\) and satisfies
\begin{equation}
 \label{eq:cmc-graph-inverse}
 \begin{aligned}
 L_i(w,d)
 &=\bigl(\Delta_{g_i}w-d,\langle w\rangle_i\bigr),\\
 L_i^{-1}(f,e)
 &=
 \left(
 e+\Delta_{g_i}^{-1}(f-\langle f\rangle_i),
 -\langle f\rangle_i
 \right),
 \end{aligned}
\end{equation}
where \(\Delta_{g_i}^{-1}\) acts on mean-zero functions.

The smooth precompactness of the metrics \(g_i\), the uniform lower
bound for their first positive scalar eigenvalues, and the uniform
Schauder estimates give
\[
 \|L_i^{-1}(f,e)\|_{\mathcal X_i}
 \le C_0\|(f,e)\|_{\mathcal Y_i},
\]
with \(C_0\) independent of \(i\). For a sufficiently small
$\rho_*>0$, independent of $i$, the derivatives on
$B_{\mathcal X_i}(0,\rho_*)$ of
\(\mathscr A_{i,a}^{\widehat G_i}\) have uniform Lipschitz bounds, and
\begin{equation}
 \label{eq:cmc-map-errors}
 \sup_{|a|\leq5}
 \left(
 \|\mathscr A_{i,a}^{\widehat G_i}(0,0)\|_{\mathcal Y_i}
 +
 \|D\mathscr A_{i,a}^{\widehat G_i}(0,0)-L_i\|_{\mathcal X_i\to\mathcal Y_i}
 \right)
 \leq C\varepsilon_i^{(3)}.
\end{equation}
The implicit function theorem therefore gives, for every \(|a|\leq4\), a unique pair $\bigl(\varphi_i(a),c_i(a)\bigr)$
in a fixed neighborhood of the origin such that
\[
 \mathscr H_i^{\widehat G_i}
 \bigl(a,\varphi_i(a)\bigr)\equiv c_i(a),
 \qquad
 \langle\varphi_i(a)\rangle_i=0.
\]
The uniqueness neighborhood is independent of \(i\) and \(a\).

The derivative of the map \eqref{eq:cmc-augmented-map} at this solution remains uniformly invertible. The implicit function theorem consequently gives, for every fixed
\(k,m\geq0\),
\begin{equation}
 \label{eq:cmc-graph-small}
 \sum_{\ell=0}^{m}
 \left(
 \|\partial_a^\ell\varphi_i\|_
 {C^{k+2,\alpha}(N_i)}
 +
 |\partial_a^\ell c_i|
 \right)
 \leq
 C_{k,m}\varepsilon_i^{(k+m+3)}
 =
 o_{R_i}(1),
 \qquad |a|\leq4.
\end{equation}
In particular, after increasing \(R_\ast\),
\[
 |\varphi_i|<\frac1{100},
 \qquad
 1+\partial_a\varphi_i>\frac12.
\]

Define
\[
 \Xi_i(a,x)
  \defeq 
 \Phi_{R_i}\bigl(a+\varphi_i(a,x),x\bigr).
\]
The last inequality shows that
\(\Xi_i:(-3,3)\times N_i\to M\) is a diffeomorphism onto its image. Put
\[
 U_i \defeq \Xi_i((-3,3)\times N_i).
\]
Using \eqref{eq:levelclose},
\eqref{eq:local-cylinder-coverage}, and
\eqref{eq:cmc-graph-small}, we obtain
\begin{equation}
 \label{eq:cmc-central-coverage}
 b^{-1}\bigl((R_i-2,R_i+2)\bigr)
 \subset U_i
 \subset
 b^{-1}\bigl((R_i-4,R_i+4)\bigr).
\end{equation}

Write
\[
 \Sigma_{i,a} \defeq \Xi_i(\{a\}\times N_i).
\]
These are the desired CMC hypersurfaces. Each $\Sigma_{i,a}$ is called a leaf.

\textbf{Claim 2: The leaves are either disjoint or entirely coincide. They are all diffeomorphic to the limit cross section for sufficiently large $R_i$.}

We show that any two intersecting leaves coincide. Suppose that \(\Sigma_{i,a}\) and \(\Sigma_{j,a'}\) meet, where \(|a|,|a'|<3\), and set
\[
 \lambda_i \defeq R_i+a,\qquad
 \lambda_j \defeq R_j+a'.
\]
\eqref{eq:levelclose} and
\eqref{eq:cmc-graph-small} imply
\begin{equation}
 \label{eq:overlap-levels}
 |\lambda_i-\lambda_j|=o(1),
 \qquad
 \sup_{\Sigma_{i,a}}|b-\lambda_j|=o(1),
\end{equation}
where the error tends to zero as
\(\min\{R_i,R_j\}\to\infty\). The same argument on the \(i\)-chart tube \(|t-a|<1\) shows that this entire tube lies in the \(j\)-th product chart. On it define
\[
 T_{ji} \defeq \Phi_{R_j}^{-1}\circ\Phi_{R_i}.
\]
The identity
\[
 T_{ji}^*\widehat G_j=\widehat G_i
\]
and uniform equivalence with the product metrics give uniform bounds
for \(DT_{ji}\) and its inverse. In coordinate patches, the
transformation rule for the Levi-Civita connections is
\[
 \partial_{\mu\nu}T_{ji}^A
 =
 (\Gamma_i)_{\mu\nu}^{\lambda}\partial_\lambda T_{ji}^A
 -
 (\Gamma_j)^A_{BC}(T_{ji})
 \partial_\mu T_{ji}^B\partial_\nu T_{ji}^C.
\]
It gives a uniform \(C^2\)-bound, and successive differentiation gives
uniform bounds in every finite order. If \(T_{ji}^0\) denotes the
radial component, the two Busemann alignment estimates yield
\[
 \|T_{ji}^0-t-(R_i-R_j)\|_{C^0((a-1,a+1)\times N_i)}=o(1).
\]
Interpolation on the smaller tube \(|t-a|<1/2\) therefore gives
\begin{equation}
 \label{eq:cmc-transition-coordinate}
 \|T_{ji}^0-t-(R_i-R_j)\|_{C^{k,\alpha}((a-1/2,a+1/2)\times N_i)}
 =o(1)
\end{equation}
for every fixed \(k\). These estimates are uniform over all
overlapping chart pairs. Otherwise a failing sequence, together with
smooth precompactness of $\{(N_j,g_j)\}$ and the exact isometry identity,
would give a contradiction.

It follows from \eqref{eq:cmc-graph-small} and
\eqref{eq:cmc-transition-coordinate} that the normal to
\(\Sigma_{i,a}\) is \(o(1)\)-close to the positive radial normal in the
\(j\)-th chart. Hence the projection
\[
 \pi_j:\Sigma_{i,a}\longrightarrow N_j
\]
is a local diffeomorphism. Since \(\Sigma_{i,a}\) is compact, its image
is both open and closed in the connected manifold \(N_j\), and
\(\pi_j\) is a finite covering.

The covering has degree one. Indeed, if $\mathrm{deg}(\pi_j)=d\ge 1$, then for each point \(x\in N_j\), order the finitely many radial $\R$-coordinates of the preimages as $t_1(x)<\ldots<t_d(x)$. The strict inequality is because the corresponding $\R$-coordinates cannot coincide as $\Sigma_{i,a}$ is embedded, and the ordering cannot change when $x$ varies. Each $t_k(x)$ defines a smooth global section of the covering by $x\mapsto (t_k(x),x)$, $x\in N_j$ and $k=1,\ldots,d$. The images of these sections are disjoint open and closed subsets of the connected hypersurface \(\Sigma_{i,a}\), so $d=1$.

Thus the entire leaf has a unique representation in the \(j\)-th
chart,
\[
 t=\widetilde a+\widetilde\varphi(x),
 \qquad
 \langle\widetilde\varphi\rangle_j=0.
\]
The estimate \eqref{eq:cmc-transition-coordinate} implies
\[
 |\widetilde a-a'|
 +
 \|\widetilde\varphi\|_{C^{2,\alpha}(N_j)}
 =o(1),
 \qquad
 |\widetilde a|<4.
\]
Since the graph has the same constant mean curvature as
\(\Sigma_{i,a}\), the pair
\((\widetilde\varphi,c_i(a))\) is a small zero of the map
\eqref{eq:cmc-augmented-map} in the \(j\)-th chart with parameter \(\widetilde a\). Uniform uniqueness
identifies it with
\((\varphi_j(\widetilde a),c_j(\widetilde a))\). The leaves in the
\(j\)-th family are disjoint, and the graph meets \(\Sigma_{j,a'}\). Hence,
\[
 \widetilde a=a',
 \qquad
 \Sigma_{i,a}=\Sigma_{j,a'}.
\]
Thus any two of the CMC leaves are either disjoint or agree as entire hypersurfaces.

\textbf{Claim 3: The CMC hypersurfaces foliate the end and can be parametrized properly by $\R$.}

 Let
\[
 U \defeq \bigcup_{i\geq0}U_i.
\]
Adjacent sets overlap by \eqref{eq:cmc-central-coverage}, so \(U\) is
connected and contains all sufficiently large Busemann level sets. Since intersecting CMC leaves agree as entire hypersurfaces, these leaves partition \(U\) into closed connected CMC hypersurfaces.

Let \(\mathscr L\) be the index space of leaves equipped with the quotient topology. A sufficiently small parameter interval in any \(U_i\) is a saturated product collar. Indeed, every leaf meeting the collar is one of its members. These collars make \(\mathscr L\) a smooth one-dimensional manifold without boundary. Distinct leaves are disjoint compact sets and have positive distance, so their saturated collars may be chosen disjoint. Hence, \(\mathscr L\) is Hausdorff. Since the quotient map \(U\to\mathscr L\) is open in these collars and \(U\) is second countable, so is \(\mathscr L\).

The index space is connected and noncompact. Indeed, a compact index space would admit a finite cover by collars lying in finitely many \(U_i\), whereas \(b\) is unbounded on \(U\). Therefore
\[
 \mathscr L\cong\mathbb R.
\]
The positive normal of $\Sigma_{i,a}=\{(a+\varphi_i(a,x),x):x\in N_i\}$, which is the normal pointing to $\{t>a+\varphi_i(a,x)\}$, defines a consistent orientation.

Every leaf with sufficiently large $\inf_\Sigma b$ can be represented as $\Sigma_{i,a}$ with $|a|<1$, after passing to an overlapping chart. 

\textbf{Claim 4: Each $\Sigma_{i,a}$ is separating.}

 Fix such a representation and put
\[
 \lambda \defeq R_i+a.
\]
Any path joining a point with \(b<\lambda-\frac12\) to a point with
\(b>\lambda+\frac12\) has a subpath contained in the Busemann tube
\[
 \lambda-\frac12\leq b\leq\lambda+\frac12.
\]
This tube lies in the \(i\)-th chart. At the two endpoints of the subpath, the signed graph function satisfies
\begin{equation}
 \label{eq:cmc-separation}
 t-a-\varphi_i(a,x)
 \leq-\frac12+o(1)<0,
 \qquad
 t-a-\varphi_i(a,x)
 \geq\frac12-o(1)>0.
\end{equation}
Hence the path meets \(\Sigma_{i,a}\). This finishes the proof of the separating property.

For each \(x\in N_i\), the curves
\[
 \alpha\longmapsto\Xi_i(\alpha,x),
 \qquad
 \alpha\in[a-1,a+1],
\]
remain in \(U_i\). Their endpoints lie on the two sides of
\(\Sigma_{i,a}\), by \eqref{eq:levelclose} the Busemann values of the endpoints are below \(\lambda-\frac12\) and above \(\lambda+\frac12\), respectively. Since the inverse images of the two components of
\(\mathscr L\setminus\{[\Sigma_{i,a}]\}\) are connected,
\eqref{eq:cmc-separation} shows that $\{b>\lambda+\tfrac12\}$ lies on the positive side of the leaf, whereas $U\cap\{b<\lambda-\tfrac12\}$ lies on its negative side.

Choose one such leaf \(\Sigma_0=\Sigma_{i,a_0}\) with associated Busemann value
\(\lambda_0>R_\ast+5\), and let \(E_+\) be the component of \(M\setminus\Sigma_0\) containing all sufficiently large Busemann level sets.
Equation~\eqref{eq:cmc-separation} gives
\begin{equation}
 \label{eq:cmc-end-separation}
 \{b>\lambda_0+\tfrac12\}
 \subset E_+
 \subset
 \{b\geq\lambda_0-\tfrac12\}.
\end{equation}
The second inclusion follows because any path from a point with \(b<\lambda_0-\frac12\) to the Busemann level set of large values must cross \(\Sigma_0\).

The position of $U_i$ \eqref{eq:cmc-central-coverage} gives $\{b>R_\ast-2\}\subset U$.
Since \(\lambda_0>R_\ast+5\), the right inclusion in \eqref{eq:cmc-end-separation} implies \(E_+\subset U\). The left  inclusion gives
\[
 M\setminus E_+
 \subset
 \{b\leq\lambda_0+\tfrac12\},
\]
which is compact by \eqref{eq:busemann-bounds}. Moreover,
\(\partial E_+=\Sigma_0\). Indeed, the inclusion \(\partial E_+\subset\Sigma_0\) is immediate, while the positive
curves \(\alpha\mapsto\Xi_i(\alpha,x)\) issuing from \(\Sigma_0\)
show the reverse inclusion.

Set
\[
 K \defeq M\setminus E_+.
\]
For sufficiently small $\varepsilon>0$, the map $(t,x)\mapsto\Xi_i(a_0+t,x)$, with $-\varepsilon<t\le0$ and $x\in N_i$, gives a collar $(-\varepsilon,0]\times\Sigma_0$ in \(K\). Thus \(K\) is a compact domain with nonempty interior. The separating property of $K$ also shows
\[
 \partial K=\Sigma_0,
 \qquad
 \overline E_+
 = E_+\cup\Sigma_0
 =
 M\setminus\operatorname{int}K.
\]

\textbf{Claim 5: The volume and normal speed bounds, and the convergence \eqref{eq:cmc-local-product-limit}, hold.}

The foliation of \(U\) restricts to \(\overline E_+\). For a leaf
\(\Sigma\subset\overline E_+\), let \(\Omega_\Sigma\) be the closed
region between \(\Sigma_0\) and \(\Sigma\), and set
\[
 \tau(\Sigma) \defeq \vol_g(\Omega_\Sigma),\qquad \tau(\Sigma_0)=0.
\]
If \(\Sigma=\Sigma_{i,a}\) and \(\lambda=R_i+a\), then
\(\Omega_\Sigma\subset\{\lambda_0-\tfrac12\leq b\leq\lambda+\tfrac12\}\),
so this region is compact. Define
\[
 \mathcal V_i \defeq (\Xi_i)_*\partial_a,
 \quad
 w_i \defeq g(\mathcal V_i,\nu_i)=1+o(1)>0,
 \quad \mathcal V_i^\top \defeq \mathcal V_i-w_i\nu_i,
\]
where \(\nu_i\) is the outward unit normal. The volume variation formula gives
\begin{equation}
 \label{eq:cmc-volume-parameter}
 \frac{d\tau}{da}=\int_{\Sigma_{i,a}}w_i\,d\vol_{\Sigma_{i,a}}>0.
\end{equation}
Thus \(\tau\) is a smooth increasing leaf coordinate and \(\tau\to\infty\).

Set \(V(\tau) \defeq \vol(\Sigma_\tau)\). After increasing \(R_\ast\),
\begin{equation}
 \label{eq:cmc-volume-bounds}
 \frac12V_\infty\leq V(\tau)\leq2V_\infty.
\end{equation}
Consequently
\[
 r(\tau) \defeq \int_0^\tau\frac{d\sigma}{V(\sigma)}
\]
is a smooth increasing coordinate with range \([0,\infty)\). In the sequel, we can and will label leaves by \(r\).

By \eqref{eq:cmc-volume-parameter},
\(dr/da=\fint_{\Sigma_{i,a}}w_i\,d\vol_{\Sigma_{i,a}}\),
so their normal speed is
\begin{equation}
 \label{eq:cmc-normal-speed}
 u=w_i\frac{da}{dr}
 =\frac{w_i}{\fint_{\Sigma_{i,a}}w_i\,d\vol_{\Sigma_{i,a}}},
 \qquad
 \frac1{V(r)}\int_{\Sigma_r}u\,d\vol_{\Sigma_r}=1.
\end{equation}
An increasing change of leaf parameter multiplies the numerator and denominator by the same factor, while changing coordinates within a leaf affects only the tangential variation. Thus \(u\) is a well-defined smooth positive function on \(\overline E_+\).

Let
\[
 \Psi:[0,\infty)\times\Sigma_0
 \longrightarrow\overline E_+
\]
be the forward flow of \(u\nu\), with \(\Psi(0,x)=x\). Along every
flow line,
\[
 \frac{d}{dr}r(\Psi(r,x))=dr(u\nu)=1,
\]
so \(\Psi(r,x)\in\Sigma_r\). The flow exists on every finite
\(r\)-interval because the corresponding trajectory remains in the
compact region swept out between \(\Sigma_0\) and \(\Sigma_R\).

For fixed \(r\), the map
\[
 \Psi_r:\Sigma_0\longrightarrow\Sigma_r
\]
is a local diffeomorphism. Its image is open, and it is closed because
\(\Sigma_0\) is compact. Since \(\Sigma_r\) is connected, the map is
surjective. Uniqueness of integral curves makes it injective, and the
backward flow gives its smooth inverse. Therefore
\[
 \Psi:[0,\infty)\times\Sigma_0
 \longrightarrow
 M\setminus\operatorname{int}K
\]
is a diffeomorphism.

By Proposition~\ref{prop:local-gauge}, \(\Sigma_0\) has the common
diffeomorphism type of the asymptotic cross sections. Fix a
diffeomorphism
\[
 \vartheta_0:N\longrightarrow\Sigma_0
\]
and set
\[
 \Phi(r,x) \defeq \Psi(r,\vartheta_0(x)).
\]
Since \(\partial_r\Phi=u\nu\) is orthogonal to the leaves,
\[
 \Phi^*g=u^2dr^2+\gamma(r),
\]
and \eqref{eq:cmc-metric} follows from
\eqref{eq:cmc-normal-speed}.

We now prove \eqref{eq:cmc-local-product-limit} by comparing the normal-transport coordinates $\Phi$ we just constructed with the product chart $\Phi_i=\Phi_{R_i}$ constructed in Proposition \ref{prop:local-gauge}. Fix a compact interval $J\subset(-3,3)$. In the $\Xi_i$ coordinates, \eqref{eq:cmc-graph-small} and
$\widehat G_i\to\overline G_i$ give, for each fixed $k$, the following
identities with errors tending to zero in $C^{k,\alpha}(J\times N_i)$:
\begin{equation}
 \label{eq:graph-product-comparison}
 \Xi_i^*g=da^2+g_i+o(1),
 \qquad
 w_i=1+o(1),
 \qquad
 \mathcal V_i^\top=o(1).
\end{equation}
It follows that
\[
 \frac{dr}{da}
 =
 \fint_{\Sigma_{i,a}}w_i
 =
 1+o(1),
 \qquad
 \frac{da}{dr}=1+o(1),
 \qquad
 u=1+o(1).
\]

At fixed leaf coordinates, the \(r\)-velocity of \(\Xi_i\) is
\[
 \frac{da}{dr}\mathcal V_i
 =
 u\nu_i+\frac{da}{dr}\mathcal V_i^\top.
\]
Let \(Z_i(r,\cdot)\) be the vector field on \(N_i\) determined by
\[
 \bigl(\Xi_i(a,\cdot)\bigr)_*Z_i
 =
 -\frac{da}{dr}\mathcal V_i^\top.
\]
Then \(Z_i=o(1)\) in every fixed norm. Let \(\phi_r\) solve
\[
 \partial_r\phi_r=Z_i(r)\circ\phi_r
\]
with \(\phi_r=\operatorname{id}\) on the central leaf, and define
\[
 \widetilde\Xi_i(r,x)
  \defeq 
 \Xi_i(a(r),\phi_r(x)).
\]
By construction,
\[
 \partial_r\widetilde\Xi_i=u\nu_i.
\]
Differentiating the ODE for \(\phi_r\) shows that \(\phi_r\) is
\(o(1)\)-close to the identity in every fixed derivative norm on a
fixed interval $r(J)$. Hence, with errors measured in
$C^{k,\alpha}(r(J)\times N_i)$,
\begin{equation}
 \label{eq:normal-product-comparison}
 \widetilde\Xi_i^*g
 =
 dr^2+g_i+o(1),
 \qquad
 u=1+o(1).
\end{equation}
The second fundamental form is \(o(1)\) as well, by the convergence \(\widehat G_i\to dt^2+g_i\), and \eqref{eq:cmc-graph-small}.

Now fix $L>0$, $k\ge0$, and $r_j\to\infty$. Choose
$x_j\in\Sigma_{r_j}$ and set $R_j=b(x_j)$.
Apply Proposition~\ref{prop:local-gauge} with a product chart of fixed width larger than $2L+4$, centered at $R_j$, and let
$(N_j,g_j)$ be its asymptotic model. The oscillation of the Busemann function on $\Sigma_{r_j}$ tends to zero by \eqref{eq:cmc-graph-small} and
\eqref{eq:levelclose}. The width larger than $2L+4$ was chosen so that the chart contains the required $r$-interval because $dr/da=1+o(1)$.
The uniform \(C^{2,\alpha}\)-uniqueness for
\eqref{eq:cmc-augmented-map}
identifies the CMC graphs in this wider chart with the restriction of the global foliation. The normal coordinate comparison above therefore provides a diffeomorphism
\[
 \vartheta_j:N_j\longrightarrow N
\]
for which
\[
 \mathcal C_j(s,x)
 =
 \Phi(r_j+s,\vartheta_j(x))
\]
agrees with the corresponding normal parametrization on \([-L,L]\times N_j\). Equations \eqref{eq:normal-product-comparison} and \eqref{eq:cmc-graph-small} then give
\eqref{eq:cmc-local-product-limit}. The higher order statements follow by applying \eqref{eq:cmc-augmented-map}.

The bounds for \(V\) were proved in \eqref{eq:cmc-volume-bounds}. On the end, \eqref{eq:cmc-local-product-limit} gives \(u\to1\). On the remaining compact region, \(u\) is smooth and strictly positive. Hence \(u\) also admits global positive upper and lower bounds. This completes the
proof.
\end{proof}

With the normalized normal transport coordinates in place, we now derive some corresponding hypersurface identities and the resulting monotonicity of the leaf volume. These formulas provide the scalar constraint and the global volume control to be used in Section~\ref{sec:uniqueness}. However, since in the next section a further coordinate change will be performed, introducing a shift vector field $X$, the actual identities being used are variants of the following ones. Nevertheless, the derivations of the identities are almost the same.

\begin{proposition}\label{prop:cmc-volume}
In the coordinates \eqref{eq:cmc-metric}, let \(\nu\) be the unit
normal pointing toward the noncompact end. Identifying
\(\partial_r\) with its pushforward under \(\Phi\), we have
\(\partial_r=u\nu\), so \(u>0\) is the normal speed of the leaf
variation. For vectors \(Y,Z\) tangent to a leaf, recall the second fundamental form and mean curvature of leaves
\[
 A(Y,Z) \defeq g(\nabla_Y\nu,Z),
 \qquad
 H \defeq \tr_\gamma A.
\]
Then \(H=H(r)\), independent of points on the leaf, and
\begin{equation}
 \label{eq:cmc-variation}
 \partial_r\gamma=2uA,\qquad
 H'=\Delta_\gamma u-u|A|_\gamma^2,\qquad
 R_\gamma=H^2-|A|_\gamma^2.
\end{equation}
Moreover,
\begin{equation}
 \label{eq:cmc-volume}
 (\log V)'=H,\qquad
 H'=-\frac1{V}\int_Nu|A|_\gamma^2\,d\vol_\gamma\leq0.
\end{equation}
Consequently,
\[
 H\geq0,\qquad H\longrightarrow0,\qquad
 V(r)\longrightarrow V_\infty,
\]
where \(V_\infty\) is the common cross-section volume in \eqref{eq:common-cross-section-volume}. Finally,
\begin{equation}
 \label{eq:cmc-volume-tail}
 \int_a^\infty
 \frac{1+r-a}{V(r)}
 \int_Nu|A|_\gamma^2\,d\vol_\gamma\,dr
 =
 H(a)+\log\frac{V_\infty}{V(a)}
 \longrightarrow0
 \qquad\text{as }a\to\infty.
\end{equation}
Thus the limiting CMC leaf volume agrees with the volume of every asymptotic cross section.
\end{proposition}

\begin{proof}
Let \(e_i\) be tangent coordinate vectors extended by the \(r\)-coordinate, so that \([\partial_r,e_i]=0\). Write
\[
 A_{ij}=g(\nabla_{e_i}\nu,e_j),
 \qquad
 \partial_r=u\nu.
\]
Torsion-freeness gives
\[
 \nabla_{\partial_r}e_i
 =\nabla_{e_i}(u\nu).
\]
Using metric compatibility and the symmetry of \(A\), we obtain
\begin{equation}
 \label{eq:cmc-metric-derivation}
 \partial_r\gamma_{ij}
 =
 g(\nabla_{e_i}(u\nu),e_j)
 +
 g(e_i,\nabla_{e_j}(u\nu))
 =
 2uA_{ij}.
\end{equation}
It follows that
\begin{equation}
 \label{eq:cmc-area-form}
 \partial_r(d\vol_\gamma)
 =
 \frac12\tr_\gamma(\partial_r\gamma)\,d\vol_\gamma
 =
 uH\,d\vol_\gamma.
\end{equation}

We next compute the variation of \(H\). Differentiating
\(g(\nu,\nu)=1\) and \(g(\nu,e_i)=0\) gives
\[
 \nabla_{\partial_r}\nu=-\nabla^\gamma u.
\]
Indeed, its normal component vanishes, while
\[
 g(\nabla_{\partial_r}\nu,e_i)
 =
 -g(\nu,\nabla_{\partial_r}e_i)
 =
 -g(\nu,\nabla_{e_i}(u\nu))
 =
 -e_i(u).
\]
With the conventions
\[
 \begin{aligned}
 R(X,Y)Z
 &=\nabla_X\nabla_YZ-\nabla_Y\nabla_XZ-\nabla_{[X,Y]}Z,\\
 \Rm(X,Y,Z,W)&=g(R(X,Y)Z,W),
 \end{aligned}
\]
differentiating \(A_{ij}\) yields
\begin{equation}
 \label{eq:cmc-second-form-variation}
 \partial_rA_{ij}
 =
 -\nabla_i^\gamma\nabla_j^\gamma u
 +uA_i{}^kA_{kj}
 -u\Rm_g(\nu,e_i,e_j,\nu).
\end{equation}
Since
\[
 \partial_r\gamma^{ij}=-2uA^{ij},
 \qquad
 \Delta_\gamma=-\tr_\gamma\nabla_\gamma^2,
\]
tracing \eqref{eq:cmc-second-form-variation} gives
\begin{align*}
 H'=
 (\gamma^{ij})'A_{ij}
 +\gamma^{ij}A'_{ij} =
 \Delta_\gamma u
 -u\bigl(|A|_\gamma^2+\Ric_g(\nu,\nu)\bigr)
 =
 \Delta_\gamma u-u|A|_\gamma^2,
\end{align*}
where the last equality uses \(\Ric_g=0\). The contracted Gauss
equation is
\[
 R_g
 =
 R_\gamma+2\Ric_g(\nu,\nu)-H^2+|A|_\gamma^2.
\]
Since \(g\) is Ricci flat,
\[
 R_\gamma=H^2-|A|_\gamma^2.
\]
This is the scalar constraint, and
\eqref{eq:cmc-variation} follows.

Because every leaf has constant mean curvature, \(H\) depends only
on \(r\). Integrating \eqref{eq:cmc-area-form} and using the
normalization in \eqref{eq:cmc-metric}, we find
\[
 V'
 =
 \int_NuH\,d\vol_\gamma
 =
 H\int_Nu\,d\vol_\gamma
 =
 HV,
\]
and hence
\[
 (\log V)'=H.
\]
Integrating the mean-curvature variation over the closed leaf gives
\[
 VH'
 =
 \int_N\Delta_\gamma u\,d\vol_\gamma
 -
 \int_Nu|A|_\gamma^2\,d\vol_\gamma
 =
 -\int_Nu|A|_\gamma^2\,d\vol_\gamma.
\]
This proves \eqref{eq:cmc-volume}.

We now use the two sided bounds for \(V\) from
Theorem~\ref{thm:cmc-foliation}. If \(H(r_0)<0\), then the monotonicity
of \(H\) gives \(H(r)\leq H(r_0)<0\) for \(r\geq r_0\), and therefore
\[
 \log V(r)
 \leq
 \log V(r_0)+H(r_0)(r-r_0)
 \longrightarrow-\infty,
\]
contrary to the positive lower bound for \(V\). Thus \(H\geq0\), and \(V\) is nondecreasing. Since \(V\) is bounded above,
\[
 V(r)\longrightarrow V_*\in(0,\infty).
\]
The nonincreasing function \(H\) also has a limit
\(H_\infty\geq0\). If \(H_\infty>0\), then
\[
 \log V(r)
 \geq
 \log V(r_0)+H_\infty(r-r_0)
 \longrightarrow+\infty,
\]
contradicting the upper bound for \(V\). Hence \(H\to0\).

To identify \(V_*\), choose any sequence \(r_j\to\infty\) and apply
\eqref{eq:cmc-local-product-limit}. The induced metrics on the central
leaves converge to \(g_j\), so
\[
 V(r_j)-\vol(N_j,g_j)\longrightarrow0.
\]
By \eqref{eq:common-cross-section-volume},
\(\vol(N_j,g_j)=V_\infty\), and therefore \(V_*=V_\infty\).
Integrating \((\log V)'=H\) now gives
\begin{equation}
 \label{eq:cmc-H-integral}
 \int_a^\infty H(r)\,dr
 =
 \log\frac{V_\infty}{V(a)}.
\end{equation}

Finally, by \eqref{eq:cmc-volume},
\[
 -H'(r)
 =
 \frac1{V(r)}
 \int_Nu|A|_\gamma^2\,d\vol_\gamma.
\]
For \(\rho>a\), integration by parts gives
\begin{equation}
 \label{eq:cmc-tail-finite}
 \int_a^\rho
 \frac{1+r-a}{V(r)}
 \int_Nu|A|_\gamma^2\,d\vol_\gamma\,dr
 =
 H(a)-(1+\rho-a)H(\rho)
 +\int_a^\rho H(r)\,dr.
\end{equation}
Since \(H\) is nonnegative, nonincreasing, and integrable by
\eqref{eq:cmc-H-integral},
\[
 \frac{\rho-a}{2}H(\rho)
 \leq
 \int_{(a+\rho)/2}^{\rho}H(r)\,dr
 \longrightarrow0.
\]
Together with \(H(\rho)\to0\), this gives
\[
 (1+\rho-a)H(\rho)\longrightarrow0.
\]
Letting \(\rho\to\infty\) in \eqref{eq:cmc-tail-finite} and using \eqref{eq:cmc-H-integral} proves \eqref{eq:cmc-volume-tail}. Its right hand side tends to zero because \(H(a)\to0\) and \(V(a)\to V_\infty\).
\end{proof}

\section{Fixed slice coordinates and shift estimates}\label{sec:fixed-slice}

Section~\ref{sec:construction-gauge} constructs the global CMC foliation and proves the identities for its logarithmic leaf volume. We now compare its leaf metrics with the prescribed integrable cross section metric \(g_N=g_0\). We first place the leaf metrics in a fixed slice satisfying a divergence condition $\delta(\gamma-g_0)=0$, then construct a Ricci flat family $g_z$ near $g_0$ with parameter $z$ so that $P_{K_0}(g_z-g_0)=z$. We then decompose $\gamma-g_z$ according to the sign of the spectrum of the Lichnerowicz Laplacian. The scalar Gauss constraint will relate the mean curvature and second fundamental form of the leaves to Lichnerowicz eigenvalues. Finally, the positive spectrum estimate in the next section will control the kernel and negative spectral components and give exponential convergence.

\subsection{Fixed slice coordinates along the foliation}

Set \(g_0 \defeq g_N\). For \(m\ge2\), let
\begin{equation}
 \label{eq:fixed-slice-space}
 \mathscr S^{m,\alpha}
  \defeq g_0+\left\{h\in C^{m,\alpha}(\operatorname{Sym}^2T^*N):  \delta_{g_0}h=0\right\}.
\end{equation}
This is an affine Banach space, with the norm inherited from
$C^{m,\alpha}(N)$. We use the balls
\begin{equation}
 \label{eq:slice-balls}
 B_{\mathscr S^{m,\alpha}}(g_0,\tau)
  \defeq \left\{\gamma\in\mathscr S^{m,\alpha}:
          \|\gamma-g_0\|_{C^{m,\alpha}(N)}<\tau\right\}.
\end{equation}
All radii below are chosen so that these tensors are positive
definite. We begin with a fixed-slice statement on a regularity
scale with spare derivatives for the pullback operation.

\begin{lemma}
\label{lem:fixed-affine-slice}
Fix \(m\geq2\) and \(\alpha\in(0,1)\), and define
\[
 \mathcal X_{\perp}^{m+1,\alpha}
  \defeq 
 \left\{
  Y\in C^{m+1,\alpha}(TN):
  Y\perp_{L^2(N)}\operatorname{Kill}(g_0)
 \right\},
\]
and
\[
 \mathcal Y_{\perp}^{m-1,\alpha}
  \defeq 
 \left\{
  \omega\in C^{m-1,\alpha}(T^*N):
  \omega\perp_{L^2(N)}
  \operatorname{Kill}(g_0)^\flat
 \right\}.
\]
There are neighborhoods \(\mathcal V_m\) of \(0\) in
\(\mathcal X_{\perp}^{m+1,\alpha}\) and \(\mathcal U_m\) of \(g_0\)
in the space of \(C^{m+2,\alpha}\) metrics such that every
\(\widetilde g\in\mathcal U_m\) admits a unique
\(Y=Y(\widetilde g)\in\mathcal V_m\) satisfying
\begin{equation}
 \label{eq:fixed-affine-slice}
 \delta_{g_0}\bigl(\varphi_Y^*\widetilde g-g_0\bigr)=0,
 \qquad
 \varphi_Y(x) \defeq \exp_x^{g_0}Y(x).
\end{equation}
Moreover,
\[
 \|Y(\widetilde g)\|_{C^{m+1,\alpha}}
 \leq
 C_m\|\widetilde g-g_0\|_{C^{m+2,\alpha}(N)},
\]
and \(\widetilde g\mapsto Y(\widetilde g)\) is \(C^1\) between the space of \(C^{m+2,\alpha}\) metrics and $\mathcal X_{\perp}^{m+1,\alpha}$.

If $\gamma$ lies in
$B_{\mathscr S^{m+1,\alpha}}(g_0,\tau_m)$,
as defined in \eqref{eq:slice-balls} with $m$ replaced by $m+1$, and $\tau_m>0$
is sufficiently small, then
\begin{equation}
 \label{eq:fixed-slice-operator}
 \mathcal B_\gamma:
 \mathcal X_{\perp}^{m+1,\alpha}
 \longrightarrow
 \mathcal Y_{\perp}^{m-1,\alpha},
 \qquad
 \mathcal B_\gamma Z \defeq \delta_{g_0}\mathcal L_Z\gamma,
\end{equation}
is an isomorphism with uniformly bounded inverse. The constructions at
different values of \(m\) are compatible. Consequently, for a smooth one-parameter family of metrics the representatives constructed here are also smooth, with estimates in every prescribed finite order.
\end{lemma}

\begin{proof}
For \(Y\) and \(\widetilde g\) in the spaces above, define
\[
 \mathfrak F_m(Y,\widetilde g)
  \defeq 
 \delta_{g_0}\bigl(\varphi_Y^*\widetilde g-g_0\bigr).
\]
For every \(K\in\operatorname{Kill}(g_0)\),
\[
 \bigl\langle
  \mathfrak F_m(Y,\widetilde g),K^\flat
 \bigr\rangle_{L^2(N)}
 =
 \bigl\langle
  \varphi_Y^*\widetilde g-g_0,
  \delta_{g_0}^*K^\flat
 \bigr\rangle_{L^2(N)}
 =0,
\]
so \(\mathfrak F_m\) takes values in
\(\mathcal Y_{\perp}^{m-1,\alpha}\).
Moreover,
\[
 \mathfrak F_m:
 \mathcal X_{\perp}^{m+1,\alpha}
 \times \mathcal U_m
 \longrightarrow
 \mathcal Y_{\perp}^{m-1,\alpha}
\]
is \(C^1\). Then at \((0,g_0)\),
\begin{equation}
 \label{eq:fixed-slice-linearization}
 D_Y\mathfrak F_m(0,g_0)[Z]
 =
 \delta_{g_0}\mathcal L_Zg_0
  =
\mathcal B_{g_0}Z
 =
 2\delta_{g_0}\delta_{g_0}^*(Z^\flat).
\end{equation}
This is a second-order elliptic, self-adjoint operator. Its kernel
consists precisely of the Killing fields, since
\[
 \left\langle
 2\delta_{g_0}\delta_{g_0}^*(Z^\flat),Z^\flat
 \right\rangle_{L^2(N)}
 =
 2\|\delta_{g_0}^*(Z^\flat)\|_{L^2(N)}^2
 =
 \frac12 \|\mathcal L_Zg_0\|^2_{L^2(N)}.
\]
It is therefore an isomorphism from
\(\mathcal X_{\perp}^{m+1,\alpha}\) onto
\(\mathcal Y_{\perp}^{m-1,\alpha}\), with the Schauder estimate
\[
 \|Z\|_{C^{m+1,\alpha}}
 \leq
 C_m
 \|\delta_{g_0}\mathcal L_Zg_0\|_{C^{m-1,\alpha}}.
\]
The Banach implicit function theorem proves
\eqref{eq:fixed-affine-slice}, uniqueness in $\mathcal V_m$, and the estimate of $\|Y(\widetilde g)\|_{C^{m+1,\alpha}}$.

For $\gamma\in
B_{\mathscr S^{m+1,\alpha}}(g_0,\tau_m)$,
the operator \(\mathcal B_\gamma\) in \eqref{eq:fixed-slice-operator} is a small perturbation
of the operator in \eqref{eq:fixed-slice-linearization}, and hence remains uniformly invertible. Uniqueness in $C^{m,\alpha}$-norm identifies the solutions obtained at higher order norms. 
\end{proof}

\begin{remark}\label{rem:fixed-slice-low-order}
There is also a formulation useful for smooth families that requires
only $C^1$-closeness of $\gamma$ to $g_0$. Put
\[
 E \defeq \mathsf{H}^1(TN)\cap\operatorname{Kill}(g_0)^\perp,\qquad
 a_\gamma(Y,Z) \defeq 
 \langle\mathcal L_Y\gamma,\delta_{g_0}^*Z^\flat\rangle_{L^2(N)}.
\]
 The dual $E^*$ is identified with the
$\mathsf{H}^{-1}$ one-forms annihilating the Killing fields.
Expanding the Lie derivative and integrating by parts on the closed
manifold $N$ gives
\[
 \tfrac12\|\mathcal L_Yg_0\|_{L^2(N)}^2
 =\|\nabla_{g_0}Y\|_{L^2(N)}^2
  +\|\operatorname{div}_{g_0}Y\|_{L^2(N)}^2.
\]
Indeed, the mixed term in the expansion satisfies
\[
 \int_N\nabla_iY_j\nabla^jY^i\,dV_{g_0}
 =\int_N(\operatorname{div}_{g_0}Y)^2\,dV_{g_0},
\]
because the Ricci term from commuting covariant derivatives vanishes.
The identity extends from smooth vector fields to
$\mathsf{H}^1(TN)$ by density. In particular, the Killing fields of
$g_0$ are precisely its parallel fields thanks to Ricci flatness.

On the $L^2$-orthogonal complement of $\operatorname{Kill}(g_0)$, the Poincar\'e inequality holds
\[
 \|Y\|_{L^2(N)}
 \le C\|\nabla_{g_0}Y\|_{L^2(N)},\qquad Y\in E.
\]
To see this, we can argue by contradiction. The failure of the inequality would give $Y_j\in E$ with $\|Y_j\|_{L^2(N)}=1$ and
$\|\nabla_{g_0}Y_j\|_{L^2(N)}\to0$.
By Rellich compactness, a subsequence converges strongly in $L^2$
and weakly in $\mathsf{H}^1$ to a vector field $Y\in E$ with $\|Y\|_{L^2(N)}=1$ and $\nabla_{g_0}Y=0$.
Elliptic regularity makes this weakly parallel field smooth, so it is a nonzero Killing field orthogonal to all Killing fields, a contradiction. Combining the Poincar\'e estimate with the preceding
identity gives
\[
 a_{g_0}(Y,Y)=\tfrac12\|\mathcal L_Yg_0\|_{L^2(N)}^2
 \ge c_0\|Y\|_{\mathsf{H}^1(N)}^2,\qquad Y\in E.
\]
Integration by parts and the formula for the Lie derivative imply
\[
 |a_\gamma(Y,Z)-a_{g_0}(Y,Z)|
 \le C\|\gamma-g_0\|_{C^1(N)}
       \|Y\|_{\mathsf{H}^1(N)}\|Z\|_{\mathsf{H}^1(N)}.
\]
Thus, if $C\|\gamma-g_0\|_{C^1(N)}\le c_0/2$, then
$a_\gamma(Y,Y)\ge(c_0/2)\|Y\|_{\mathsf{H}^1(N)}^2$.
The Lax--Milgram theorem gives a unique
$Y\in E$  for every $f\in E^*$ satisfying $\mathcal B_\gamma Y=f$, with
\[
 \|Y\|_{\mathsf{H}^1(N)}
 \le \frac{2}{c_0}\|f\|_{\mathsf{H}^{-1}(N)}.
\]
 For smooth $\gamma$,
elliptic regularity then gives an isomorphism from the fixed
$C^{\ell+2,\alpha}$ complement onto the fixed $C^{\ell,\alpha}$
complement for every $\ell\ge0$. Its higher order bounds may depend
on higher norms of $\gamma$.

The implicit function argument can then be based at any smooth $\gamma$ in the fixed slice instead of $g_0$ and in this sufficiently small $C^1$-neighborhood of $g_0$, with $\mathfrak F_m(0,\gamma)=0$ and derivative $\mathcal B_\gamma$. It gives rise to representatives for nearby smooth metrics without requiring $\gamma-g_0$ to be small in higher $C^{m+1,\alpha}$-norm. Elliptic bootstrapping of the slice equation shows that each representative is smooth.
\end{remark}

We now introduce a tensor decomposition at \(g_0\). Put
\[
 K_0 \defeq \ker\Delta^L_{g_0},
 \qquad
 S \defeq \ker\delta_{g_0},
 \qquad
 W \defeq S\cap K_0^\perp,
\]
where orthogonality is taken in \(L^2(N)\). On $K_0$ set
$|z| \defeq \|z\|_{L^2(N)}$ and define
\begin{equation}
 \label{eq:kernel-balls}
 B_{K_0}(0,\rho) \defeq \{z\in K_0:|z|<\rho\}.
\end{equation}
Denote by
\(P_{K_0}\), \(P_S\), and \(P_W\) the corresponding
\(L^2(N)\)-orthogonal projections.

The commutation identities $\delta_{g_0}\Delta^L_{g_0}=\Delta_{H,g_0}\delta_{g_0}$, $\tr_{g_0}\Delta^L_{g_0}=\Delta_{g_0}\tr_{g_0}$ imply that \(K_0\subset S\). Indeed, if \(k\in K_0\), then
\(\tr_{g_0}k\) is constant and \(\delta_{g_0}k\) is a harmonic
one-form. Since \(g_0\) is Ricci flat, \(\delta_{g_0}k\) is parallel.
For every parallel one-form \(\omega\),
\[
 \langle\delta_{g_0}k,\omega\rangle_{L^2(N)}
 =
 \langle k,\delta_{g_0}^*\omega\rangle_{L^2(N)}
 =0.
\]
Taking \(\omega=\delta_{g_0}k\) gives
\(\delta_{g_0}k=0\).

It follows that \(\Delta^L_{g_0}\) preserves \(S\) and \(W\), and
its restriction to \(W\) is invertible. The projections are bounded
on every Hölder and Sobolev space used below. More explicitly,
\begin{equation}
 \label{eq:slice-projections}
 P_ST
 =
 T-\delta_{g_0}^*
   (\delta_{g_0}\delta_{g_0}^*)^{-1}\delta_{g_0}T,
 \qquad
 P_W=P_S-P_{K_0},
\end{equation}
where the inverse is taken on the orthogonal complement of the Killing
one-forms. The asserted bounds follow from elliptic estimates, while
\(P_{K_0}\) has finite rank.

Decompose
\[
 W=W_-\oplus W_+
\]
into the negative and positive spectral subspaces of
\(\Delta^L_{g_0}\), with projections \(P_-\) and \(P_+\).
The negative subspace is finite dimensional, and \(P_+=P_W-P_-\) has the same Hölder and Sobolev boundedness as
\(P_W\). The corresponding spectra of $\Delta^L_{g_0}$ restricted to \(W_-\) and \(W_+\) are separated from zero. Moreover, every negative eigentensor is trace-free. Indeed, if
\[
 \Delta^L_{g_0}k=\lambda k,
 \qquad
 \lambda<0,
\]
then
\[
 \Delta_{g_0}(\tr_{g_0}k)
 =
 \lambda\,\tr_{g_0}k,
\]
and the nonnegativity of the scalar Laplacian forces \(\tr_{g_0}k=0\). In what follows, we write
\[
 h_\pm \defeq P_\pm h.
\]
All these spaces and projections are fixed at \(g_0\).

We next show that the integrability in
Definition~\ref{def:int} produces a finite dimensional family of Ricci flat
metrics in the fixed slice, parametrized by a ball in \(K_0\).

\begin{lemma}
\label{lem:ricci-flat-decomposition}
Fix \(m\ge4\). There are \(\varepsilon>0\) and an analytic map
\[
 \eta:\kernelball{\varepsilon}\longrightarrow W\cap C^{m,\alpha},
 \qquad \eta(0)=0,\qquad D\eta(0)=0,
\]
such that \(g_z \defeq g_0+z+\eta(z)\) is a smooth Ricci-flat metric for
every \(z\in\kernelball{\varepsilon}\). The family depends smoothly
on \(z\) in every higher norm on a common smaller parameter ball,
is equivariant under \(\operatorname{Isom}(g_0)\), and satisfies
\begin{equation}
 \label{eq:slice-family}
 \delta_{g_0}(g_z-g_0)=0,
 \qquad
 P_{K_0}(g_z-g_0)=z.
\end{equation}
Every \(\gamma\in\mathscr S^{m,\alpha}\) sufficiently close to \(g_0\)
has the unique decomposition
\[
 \gamma=g_z+h,\qquad
 z=P_{K_0}(\gamma-g_0),\qquad
 h\in W\cap C^{m,\alpha}.
\]
Moreover, \(\gamma\) is Ricci-flat if and only if \(h=0\), so the
family contains every sufficiently nearby Ricci-flat metric in the
fixed slice.
\end{lemma}

\begin{proof}
Consider
\[
 \mathcal P(z,w)
  \defeq 
 P_W\Ric(g_0+z+w),
 \qquad
 z\in K_0,\quad
 w\in W\cap C^{m,\alpha}.
\]
By \eqref{eq:deform},
\[
 D_w\mathcal P(0,0)
 =
 \frac12\Delta^L_{g_0}:
 W\cap C^{m,\alpha}
 \longrightarrow
 W\cap C^{m-2,\alpha}.
\]
This is an isomorphism. The signs of its nonzero eigenvalues are
irrelevant. The Hessian term in the linearization of Ricci curvature disappears after projection because $\nabla^2f=\delta_{g_0}^*df$ is \(L^2(N)\)-orthogonal to $S$.
The analytic implicit function theorem therefore gives an analytic map
\[
 \eta:
 K_0\supset\kernelball{\varepsilon}
 \longrightarrow W\cap C^{m,\alpha}
\]
such that
\[
 \eta(0)=0,
 \quad
 D\eta(0)=0,
 \quad
 P_W\Ric(g_0+z+\eta(z))=0.
\]

Let \(v\in K_0\). By integrability, there is a smooth Ricci flat curve
\( g_v(t)\) with
\[
 g_v(0)=g_0,
 \qquad
 \left.\frac{d}{dt}\right|_{t=0} g_v(t)=v.
\]
Apply Lemma~\ref{lem:fixed-affine-slice} to place this curve in the fixed slice:
\[
 \widetilde g_v(t)
  \defeq 
 \varphi_{Y_v(t)}^* g_v(t),
 \qquad
 \delta_{g_0}\bigl(\widetilde g_v(t)-g_0\bigr)=0.
\]
Its initial velocity remains \(v\). Indeed,
\[
 \dot{\widetilde g}_v(0)-v
 =
 \mathcal L_{\dot Y_v(0)}g_0
 \in S^\perp,
\]
whereas both \(\dot{\widetilde g}_v(0)\) and \(v\) lie in \(S\).
Hence, $\dot{\widetilde g}_v(0)=v$.

Write
\[
 \widetilde g_v(t)
 =
 g_0+z_v(t)+w_v(t),
 \qquad
 z_v(t) \defeq P_{K_0}\bigl(\widetilde g_v(t)-g_0\bigr),
 \qquad
 w_v(t)\in W.
\]
Then
\[
 z_v(t)=tv+o(t).
\]
Since \(\widetilde g_v(t)\) is Ricci flat, it satisfies
\(\mathcal P(z_v(t),w_v(t))=0\). Uniqueness in the implicit function theorem gives $w_v(t)=\eta(z_v(t))$.

Define the analytic map
\[
 \mathcal R(z) \defeq \Ric(g_0+z+\eta(z)).
\]
For every \(v\in K_0\), the preceding curve satisfies
\(\mathcal R(z_v(t))=0\) and \(z_v(t)=tv+o(t)\).
If \(\mathcal R\) were not identically zero near \(0\), let
\(\mathcal R_d\) be its first nonzero Taylor polynomial. $\mathcal R_d$ is a homogeneous polynomial of degree $d$.
Then
\[
 0=\mathcal R(z_v(t))
   =t^d\mathcal R_d(v)+o(t^d)
\]
for every \(v\), forcing \(\mathcal R_d=0\), a contradiction.
Thus \(g_z \defeq g_0+z+\eta(z)\) is Ricci flat for every sufficiently small \(z\). The slice and projection identities in the statement follow from \(z\in K_0\subset S\) and \(\eta(z)\in W\).
Every nearby Ricci-flat metric in the fixed slice can be written as \(g_0+z+w\), with \(z\in K_0\) and \(w\in W\), and satisfies
\(\mathcal P(z,w)=0\). The uniqueness in the implicit function theorem therefore gives \(w=\eta(z)\), proving that this metric belongs to the family.

Uniqueness makes the family equivariant under $\operatorname{Isom}(g_0)$. To justify regularity on a common parameter neighborhood, each $g_z$ solves the system
\[
 \operatorname{Ric}(g)+\delta_{g_0}^*\delta_{g_0}(g-g_0)=0.
\]
By \eqref{eq:deform}, the linearization of $2\operatorname{Ric}(g)+2\delta_{g_0}^*\delta_{g_0}(g-g_0)$ at $g_0$ is $\Delta^L_{g_0}h-\nabla^2\operatorname{tr}_{g_0}h$.
The principal symbol of the linearized operator at a nonzero covector $\xi$ is
\[
 h\longmapsto |\xi|^2h+\xi\otimes\xi\,\operatorname{tr}_{g_0}h .
\]
Taking the trace shows that its kernel has zero trace, and then the
symbol equation gives $h=0$. Thus the system is elliptic, and remains elliptic for metrics close to $g_0$. Elliptic estimate gives smoothness of $g_z$. differentiating the same system in $z$ gives smooth dependence in every finite order on a common smaller subset of $\kernelball{\varepsilon}$.

Finally, for a metric \(\gamma\) in the stated neighborhood, set
\(z=P_{K_0}(\gamma-g_0)\) and \(h=\gamma-g_z\).
The slice and projection identities give
\(\delta_{g_0}h=0\) and \(P_{K_0}h=0\), so \(h\in W\).
Conversely, projecting any such decomposition onto \(K_0\)
determines \(z\), and hence \(h\). The characterization of the
Ricci-flat metrics follows from the completeness of the family
proved above.
\end{proof}

We specify the coordinate terminology before applying the fixed-slice construction. Write the normal-transport map from Theorem~\ref{thm:cmc-foliation} as $\widehat\Phi$, with $\widehat\Phi^*g=\widehat u^{\,2}dr^2+\widehat\gamma(r)$.

For an interval $I$, choose for each $r\in I$ a diffeomorphism $\phi_r:N\to N$, depending smoothly on $(r,x)$. Define
\begin{equation}
 \label{eq:leafwise-pullback}
 \begin{gathered}
 F_\phi(r,x) \defeq (r,\phi_r(x)),\qquad
 \Phi_\phi \defeq \widehat\Phi\circ F_\phi,\\
 \gamma(r) \defeq \phi_r^*\widehat\gamma(r),\qquad
 u(r,x) \defeq \widehat u(r,\phi_r(x)),\\
 X(r,x) \defeq (D\phi_r(x))^{-1}\partial_r\phi_r(x).
 \end{gathered}
\end{equation}
A \emph{slice lift} of the leaf metric path on $I$ means such a
smooth choice $\phi_r$ for which
$\gamma(r)\in\mathscr S^{m,\alpha}$ lies in the chosen coordinate
neighborhood of Lemma~\ref{lem:ricci-flat-decomposition}. We include the initial
identification $\phi_a$ when specifying a lift.
The tangential vector field $X$ in \eqref{eq:leafwise-pullback} is
its \emph{shift}. Directly differentiating $F_\phi$ gives
\[
 \Phi_\phi^*g
 =u^2dr^2+\gamma_{ij}(dx^i+X^i\,dr)(dx^j+X^j\,dr),
 \qquad \partial_r=u\nu+X .
\]
A lift is \emph{normalized} if
$X(r)\perp_{L^2(N)}\operatorname{Kill}(g_0)$ for every $r$.
The original normal-transport coordinates have no shift, i.e. $X=0$. The
$r$-dependence of the leafwise diffeomorphisms generally introduces a nonzero shift while leaving the CMC leaves unchanged.

We now apply this framework to the CMC leaf metrics.

\begin{proposition}\label{prop:slice-path}
Fix $m\ge6$ and $\sigma>0$ so small that $\sliceball{2\sigma}$
is contained in the fixed slice coordinate neighborhood.
There exists $a_0$ with the following property.
\begin{enumerate}
\item\label{item:existence-of-slice} 
Starting from a leaf $\Sigma_a$, $a\ge a_0$, with an initial metric representative in $\sliceball{\sigma}$, there is an $\epsilon>0$ and an interval $[a,a+\varepsilon)$ on which one can choose a diffeomorphism $\phi_r:N\to N$ for each leaf $\Sigma_r$, smoothly in $(r,x)$, such that
\[
 \Phi_\phi(r,x)=\widehat\Phi(r,\phi_r(x)),\qquad
 \gamma(r)=\phi_r^*\widehat\gamma(r)\in\sliceball{2\sigma}.
\]
This family is a slice lift of the leaf metric path.
For as long as the representatives lie in $\sliceball{2\sigma}$, the pulled-back metric $G=\Phi_\phi^*g$ takes the form
\begin{equation}
 \label{eq:slice-metric}
 \begin{gathered}
 G=
 u^2dr^2+
 \gamma_{ij}(dx^i+X^i\,dr)(dx^j+X^j\,dr),\\
 \gamma=g_z+h, \qquad  h\in W=\ker\delta_{g_0}\cap K_0^\perp .
 \end{gathered}
\end{equation}
The leaf diffeomorphisms may be chosen with the normalization so that $X(r)\perp_{L^2(N)}\operatorname{Kill}(g_0)$.

\item\label{item: extension-of-slice} If a normalized lift is defined on $[a,T)$, $T<\infty$, and
$\gamma(r)\in\overline{\sliceball{\sigma}}$ for $a\le r<T$,
then it extends smoothly to $[a,T+\varepsilon_T)$ for some
$\varepsilon_T>0$, agreeing with the original lift on $[a,T)$.

Moreover, fix \(L>0\) and \(k\geq0\). Suppose that \(r_j\to\infty\), and that $J_j \defeq [r_j-L,r_j+L]$ is contained in the interval on which the normalized lift is defined, and that the metric representatives on $J_j$ belong to $\overline{\sliceball{\sigma}}$. Then
\begin{equation}
 \label{eq:slice-block-smallness}
 \|h\|_{C^{k,\alpha}(J_j\times N)}
 +\|\partial_r\gamma\|_{C^{k,\alpha}(J_j\times N)}
 +\|u-1\|_{C^{k,\alpha}(J_j\times N)}
 +\|X\|_{C^{k,\alpha}(J_j\times N)}
 \longrightarrow0 .
\end{equation}
 The convergence is uniform over all such intervals and any normalized lift, independent of the total length of the lift interval.
\end{enumerate}

\end{proposition}

\begin{proof}
We first prove \eqref{item:existence-of-slice}. Apply the implicit function construction of Remark~\ref{rem:fixed-slice-low-order}, based at the initial smooth slice representative, to the smooth path of leaf metrics on a short interval. If \(\gamma\) denotes the resulting representative, set
\[
 z \defeq P_{K_0}(\gamma-g_0),
 \qquad
 h \defeq \gamma-g_z.
\]
Lemma~\ref{lem:ricci-flat-decomposition} gives
\[
 \delta_{g_0}h=0,
 \qquad
 P_{K_0}h=0,
\]
and hence \eqref{eq:slice-metric}.

It remains to impose the orthogonal normalization condition.
Pullback by a \(g_0\)-isometry preserves the fixed slice, the spaces
\(K_0\), \(W\), and the family \(\{g_z\}\), and hence preserves the
decomposition in \eqref{eq:slice-metric}. Let \(\chi(r)\) solve
\[
 \partial_r\chi(r)
 =-\bigl(P_{\operatorname{Kill}(g_0)}X(r)\bigr)\circ\chi(r),
 \qquad \chi(a)=\operatorname{id}.
\]
This is an ODE on the compact group \(\operatorname{Isom}(g_0)\).
Replacing \(\phi_r\) by \(\phi_r\circ\chi(r)\) gives the shift
\[
 \chi(r)^*X+(D\chi(r))^{-1}\partial_r\chi(r)
 =\chi(r)^*\bigl(X-P_{\operatorname{Kill}(g_0)}X\bigr),
\]
which is orthogonal to \(\operatorname{Kill}(g_0)\). This proves the
local existence of a normalized lift.

Now we turn to the proof of \eqref{item: extension-of-slice}. Use the fixed $m\ge6$ and $\sigma$ from the statement, so that
\[
 \overline{\sliceball{\sigma}}\subset\sliceball{2\sigma}.
\]
This inclusion does not assert that
$\overline{\sliceball{\sigma}}$ is compact. To prove the convergence assertion, take \(r_j\to\infty\) and normalized lifts on intervals containing \(J_j=[r_j-L,r_j+L]\), with representatives in \(\overline{\sliceball{\sigma}}\) on \(J_j\). Set \(s=r-r_j\). This reparametrization is referred to as recentering (at $r_j$). By \eqref{eq:cmc-local-product-limit} and the smooth precompactness in Proposition~\ref{prop:local-gauge}, after passing to a subsequence and identifying the relevant cross sections with \(N\),
the normal transport coordinates from Theorem \ref{thm:cmc-foliation} on
\([-2L,2L]\times N\) satisfy
\begin{equation}
 \label{eq:slice-recentered-cmc}
 \begin{gathered}
 \widehat G_j
 =
 \widehat u_j^{\,2}ds^2+\widehat\gamma_j(s),\\
 \|\widehat\gamma_j-g_\infty\|_
 {C^{\ell,\alpha}([-2L,2L]\times N)}
 +
 \|\widehat u_j-1\|_
 {C^{\ell,\alpha}([-2L,2L]\times N)}
 \longrightarrow0,\\
 \|\widehat A_j\|_
 {C^{\ell,\alpha}([-2L,2L]\times N)}
 \longrightarrow0
 \end{gathered}
\end{equation}
for every fixed \(\ell\geq0\), where \(g_\infty\) is Ricci flat and
independent of \(s\).

At the central leaves there are diffeomorphisms
\(\psi_j:N\to N\) such that
\[
 \gamma(r_j)=\psi_j^*\widehat\gamma_j(0).
\]
Since the metrics on both sides are uniformly equivalent to $g_0$ and uniformly bounded in $C^{m,\alpha}$-norm, this isometry identity gives uniform bounds for \(D\psi_j\), \(D\psi_j^{-1}\). In local coordinates,
\[
 \partial_b\partial_c\psi_j^a
 =
 \Gamma^d_{bc}\bigl(\gamma(r_j)\bigr)\partial_d\psi_j^a
 -
 \Gamma^a_{ef}\bigl(\widehat\gamma_j(0)\bigr)
 \partial_b\psi_j^e\partial_c\psi_j^f.
\]
Differentiating this identity gives uniform
\(C^{m+1,\alpha}\)-bounds for \(\psi_j\) and their inverses.

Higher order bounds follow from the fixed slice equation
\begin{equation}
 \label{eq:slice-central-elliptic}
 \delta_{g_0}
 \bigl(\psi_j^*\widehat\gamma_j(0)-g_0\bigr)=\delta_{g_0}\bigl(\gamma(r_j)-g_0\bigr)=0.
\end{equation}
This equation is provided by the assumption. This is a second-order quasilinear elliptic system for \(\psi_j\).
For a variation represented by a vector field \(Y\),
its principal term in its linearization is
\[
 Y\longmapsto
 \delta_{g_0}\mathcal L_Y\gamma(r_j),
\]
which is a uniformly strongly elliptic perturbation of \(Y\mapsto\delta_{g_0}\mathcal L_Yg_0\) on the complement of \(\operatorname{Kill}(g_0)\). The bounds up to $C^{m,\alpha}$ norm just obtained, the smooth bounds for \(\widehat\gamma_j(0)\), and Schauder
estimates give uniform bounds for \(\psi_j\) in every $C^{\ell,\alpha}$ for fixed $\ell$. After passing to a further subsequence,
\begin{equation}
 \label{eq:slice-central-limit}
 \psi_j\longrightarrow\psi_\infty,
 \qquad
 \gamma(r_j)\longrightarrow
 g_* \defeq \psi_\infty^*g_\infty
\end{equation}
smoothly on \(N\). The inverse bounds imply that \(\psi_\infty\) is a diffeomorphism. Since \(g_*\) is Ricci flat and lies in the fixed slice, Lemma~\ref{lem:ricci-flat-decomposition} gives \(g_*=g_{z_*}\), where \(z_* \defeq P_{K_0}(g_*-g_0)\).
Thus the smooth compactness of representatives $\gamma(r_j)$ follows
from \eqref{eq:slice-central-elliptic}.

Use \(\psi_j\), independently of \(s\), on the whole interval $[-L,L]$ and set
\[
\overline\gamma_j(s) \defeq \psi_j^*\widehat\gamma_j(s).
\]
Although \(\overline\gamma_j(0)=\gamma(r_j)\) lies in the fixed slice, \(\overline\gamma_j(s)\) need not remain in the slice for \(s\ne0\). However \eqref{eq:slice-recentered-cmc} and
\eqref{eq:slice-central-limit} imply that they converge to the same limit as $\gamma(r_j)$ as $j\to\infty$, i.e.
\[
 \overline\gamma_j\longrightarrow g_*
\]
smoothly on \([-L,L]\times N\). Applying
Lemma~\ref{lem:fixed-affine-slice} and Remark~\ref{rem:fixed-slice-low-order} to $\overline\gamma_j(s)$ gives diffeomorphisms \(f_j(s)\), with \(f_j(0)=\operatorname{id}\), for
which \(\delta_{g_0}(f_j(s)^*\overline\gamma_j(s))=0\). The construction in Remark~\ref{rem:fixed-slice-low-order}, based at \(g_*\), gives $f_j\longrightarrow\operatorname{id}$ smoothly on \([-L,L]\times N\).

Canceling the Killing field components by composing \(g_0\)-isometries $\chi_j(s)$ with $\chi_j(0)=\operatorname{id}$ imposes
\[
 X(s)\perp_{L^2(N)}\operatorname{Kill}(g_0).
\]
The Killing vector fields generating the isometries $\chi_j(s)$ tend to zero, so $\chi_j$ itself converges smoothly to the identity on $[-L,L]$. Denote the fields obtained after the composition of $\chi_j$ by $\widetilde\gamma_j$, $\widetilde u_j$, $\widetilde X_j$, and set
\[
 \widetilde z_j
  \defeq 
 P_{K_0}(\widetilde\gamma_j-g_0),
 \qquad
 \widetilde h_j
  \defeq 
 \widetilde\gamma_j-g_{\widetilde z_j}.
\]
For brevity write
\[
 \mathcal C_L^{k,\alpha}
  \defeq 
 C^{k,\alpha}([-L,L]\times N).
\]
The convergence of the recentered fields in
\eqref{eq:slice-recentered-cmc} and of the isometries $\chi_j$, together with Lemma~\ref{lem:ricci-flat-decomposition}, gives, for every fixed \(k\),
\begin{equation}
 \label{eq:slice-fixed-interval-limit}
 \begin{aligned}
 &\|\widetilde\gamma_j-g_{z_*}\|_{\mathcal C_L^{k+1,\alpha}}
 +\|\widetilde u_j-1\|_{\mathcal C_L^{k,\alpha}}
 +\|\widetilde X_j\|_{\mathcal C_L^{k,\alpha}}
 \longrightarrow0,\\
 &\|\widetilde z_j-z_*\|_{C^{k+1,\alpha}([-L,L])}
 +\|\widetilde h_j\|_{\mathcal C_L^{k,\alpha}}
 +\|\partial_s\widetilde\gamma_j\|_{\mathcal C_L^{k,\alpha}}
 \longrightarrow0.
 \end{aligned}
\end{equation}

We now identify this auxiliary normalized lift, in particular the fields with tildes, with the given normalized lift, using their common identification at $\gamma(r_j)$.
The local uniqueness in Lemma~\ref{lem:fixed-affine-slice} and Remark~\ref{rem:fixed-slice-low-order} implies that their transition diffeomorphism lies in \(\operatorname{Isom}(g_0)\) near that leaf. Indeed, a transition near the identity factors as \(\varphi_Y\circ\chi\), with \(Y\perp\operatorname{Kill}(g_0)\) and \(\chi\in\operatorname{Isom}(g_0)\). Invariance of the slice under \(\chi\) and uniqueness of the small correction based at the given slice metric force \(Y=0\).
The shifts consequently satisfy
\[
 \widetilde X=\chi^*X+Y_\chi,
 \qquad
 Y_\chi \defeq (D\chi)^{-1}\partial_r\chi
 \in\operatorname{Kill}(g_0).
\]
Pullback by \(\chi\) preserves the Killing-field complement. Since
both shifts lie in that complement, \(Y_\chi=0\). The common central
identification gives \(\chi(0)=\operatorname{id}\), so the two lifts
agree. Applying this local uniqueness successively identifies them
throughout their common interval.

It follows that \eqref{eq:slice-fixed-interval-limit} and the identification prove
\eqref{eq:slice-block-smallness} for the given normalized lift. If this convergence were not uniform over the stated intervals and normalized lifts, a violating sequence would admit the convergent subsequence just constructed, giving a contradiction.

Suppose now that a normalized lift is defined on \([a,T)\), with \(T<\infty\), and that its representatives remain in
\(\overline{\sliceball{\sigma}}\). The normal-transport coordinates of Theorem~\ref{thm:cmc-foliation} are smooth across \(T\). Recenter them at leaves \(r_j\uparrow T\). The elliptic estimates
for \eqref{eq:slice-central-elliptic} give uniform bounds for the \(\psi_j\) in every fixed order, making the corresponding representatives \(\gamma(r_j)=\psi_j^*\widehat\gamma_j(0)\) smoothly precompact.
The closure of the representatives near \(T\) is therefore covered by finitely many fixed slice neighborhoods by Lemma~\ref{lem:fixed-affine-slice} and
Remark~\ref{rem:fixed-slice-low-order}.
The uniform bounds for \(\psi_j\) give uniform derivative bounds for
\(\psi_j^*\widehat\gamma_j(s)\). Together with the finite cover, these bounds provide a common positive interval for the local
slice constructions centered at leaves approaching \(T\). Choosing a leaf closer to \(T\) than this interval length and using the uniqueness just proved extends the given normalized lift to \([a,T+\varepsilon_T)\) for some \(\varepsilon_T>0\). This proves the continuation assertion.
\end{proof}

The next corollary is a direct consequence of \eqref{eq:slice-block-smallness} and will be used in the sequel.

\begin{corollary}
    Fix \(L>0\) and \(\varepsilon>0\). After choosing \(\rho>0\) sufficiently small and \(a\) sufficiently large, we have
\begin{equation}
 \label{eq:slice-coefficient-smallness}
 \|\gamma-g_0\|_{C^4(I\times N)}
 +\|u-1\|_{C^4(I\times N)}
 +\|X\|_{C^4(I\times N)}
 \leq\varepsilon
\end{equation}
for every interval \(I\subset[a,\infty)\) of length \(L\) contained in a normalized lift with
\(\gamma(r)\in\overline{\sliceball{\sigma}}\) and \(|z(r)|\leq\rho\) on \(I\). 
\end{corollary}

\begin{proof}
For derivatives of \(\gamma\) involving \(r\), use \eqref{eq:slice-block-smallness} with \(k=3,4\). For other derivatives along the cross section, use \(\gamma=g_z+h\) and choose \(\rho\) to be sufficiently small.
\end{proof}

\begin{lemma}\label{lem:quantitative-shift}
Let $(\gamma,u,X)$ be a normalized slice lift from
Proposition~\ref{prop:slice-path} on an interval $I$, with
$\gamma(r)\in\overline{\sliceball{\sigma}}$, and let $A$ be the
second fundamental form of its CMC leaves. Then
\begin{equation}
 \label{eq:slice-shift-equation}
 \mathcal B_\gamma X=-2\delta_{g_0}(uA).
\end{equation}
Moreover,
\begin{equation}\label{eq:slice-shift-low-order}
 \|X\|_{\mathsf{H}^1(N)}
 \le C\|\delta_{g_0}(uA)\|_{\mathsf{H}^{-1}(N)} \le C\|uA\|_{L^2(N)},
\end{equation}
and, for every integer $\ell\ge0$,
\begin{equation}\label{eq:slice-shift-schauder}
 \|X(r)\|_{C^{\ell+2,\alpha}(N)}
 \le C_\ell\|uA(r)\|_{C^{\ell+1,\alpha}(N)}.
\end{equation}
For every integer $p\ge1$, the mixed derivatives satisfy
\begin{align*}
 \|\partial_r^pX\|_{C^{\ell+2,\alpha}(N)}
 \le C_{\ell,p}\left(
   \|\partial_r^p(uA)\|_{C^{\ell+1,\alpha}(N)}+\sum_{j=1}^p
   \|\partial_r^j\gamma\|_{C^{\ell+2,\alpha}(N)}
   \|\partial_r^{p-j}X\|_{C^{\ell+2,\alpha}(N)}\right).
\end{align*}
The constant $C$ depends only on the fixed background and its $C^1$ neighborhood. The constants $C_\ell$ and $C_{\ell,p}$ also depend on the corresponding uniform coefficient bounds supplied by Proposition~\ref{prop:slice-path}, and are
independent of the total length of $I$.
\end{lemma}

\begin{proof}
Recall $\mathcal B_\gamma Y=\delta_{g_0}\mathcal L_Y\gamma$ from
\eqref{eq:fixed-slice-operator}. Differentiating $\delta_{g_0}(\gamma-g_0)=0$ and using the variation $ \gamma'=2uA+\mathcal L_X\gamma$ (see \eqref{eq:slice-metric}) gives
\eqref{eq:slice-shift-equation}. 
Since $X$ lies in the fixed Killing field complement, Remark~\ref{rem:fixed-slice-low-order} gives the first inequality in \eqref{eq:slice-shift-low-order}. The last inequality follows from the dual bound
\[
 |\langle\delta_{g_0}(uA),\omega\rangle|
 \le C\|uA\|_{L^2(N)}\|\omega\|_{\mathsf{H}^1(N)}.
\]

The elliptic estimates for \eqref{eq:slice-central-elliptic} and the
convergence \eqref{eq:slice-block-smallness} give uniform coefficient bounds for \(\mathcal B_\gamma\) in every fixed order while the representatives remain in \(\overline{\sliceball{\sigma}}\).
The Schauder estimate gives
\begin{align*}
 \|X(r)\|_{C^{\ell+2,\alpha}(N)}
 &\le C_\ell\left(
       \|\mathcal B_\gamma X(r)\|_{C^{\ell,\alpha}(N)}
       +\|X(r)\|_{L^2(N)}\right)\\
 &\le C_\ell\left(
       \|\delta_{g_0}(uA(r))\|_{C^{\ell,\alpha}(N)}
       +\|uA(r)\|_{L^2(N)}\right)\\
 &\le C_\ell\|uA(r)\|_{C^{\ell+1,\alpha}(N)} .
\end{align*}

For the mixed derivatives, differentiating the first identity in
\eqref{eq:slice-shift-equation} $p$ times gives
\begin{equation}
 \label{eq:slice-shift-r-derivatives}
 \mathcal B_\gamma\partial_r^pX
 =-2\delta_{g_0}\partial_r^p(uA)
  -\sum_{j=1}^p\binom pj
       \mathcal B_{\partial_r^j\gamma}\partial_r^{p-j}X .
\end{equation}
Every $\partial_r^pX$ remains orthogonal to the same Killing fields. Applying the inverse estimate to this identity and
using the product rule proves the stated mixed derivative bound.
Together with the uniform smooth bounds already proved, this
justifies the mixed derivative estimates by induction. Their
convergence to zero also follows from
\eqref{eq:slice-block-smallness}.
\end{proof}

The proof of Theorem~\ref{thm:main} will use these estimates to exclude
finite exits from the slice and kernel neighborhoods.

\section{Uniqueness and exponential convergence}\label{sec:uniqueness}

\subsection{The Gauss constraint and Einstein--Hilbert functional estimates}

In the coordinates \eqref{eq:slice-metric}, we introduced a shift by the vector field $X$ to impose a divergence condition. Now the second fundamental form becomes $A=(\gamma'-\mathcal L_X\gamma)/(2u)$ by \eqref{eq:leafwise-pullback} and \eqref{eq:cmc-metric-derivation}.
The Gauss and contracted Codazzi equations, together with the mean curvature variation in Proposition~\ref{prop:cmc-volume}, give, in the new coordinates,
\begin{equation}
 \label{eq:volume-evolution}
 \begin{aligned}
 \gamma'&=2uA+\mathcal L_X\gamma,&\qquad \delta_\gamma A&=0,\\
 H'&=\Delta_\gamma u-u|A|_\gamma^2,&
 R_\gamma&=H^2-|A|_\gamma^2.
 \end{aligned}
\end{equation}
Here $H$ is still constant on each leaf and the normal speed still satisfies the same normalization $\fint_N u\,dV_\gamma=1$.
In particular, the contracted Codazzi identity gives $\delta_\gamma A=-dH-\operatorname{Ric}_g(\nu,\cdot)|_{TN}=0$.
Both terms vanish because the leaves are CMC and the ambient metric is Ricci flat.
The norms on each leaf are uniformly equivalent for $\gamma\in\sliceball{2\sigma}$ and sufficiently small $\sigma$.

We will derive a relation between the second fundamental form $A$, mean curvature $H$, the kernel term $z'$ and the transverse terms $h_\pm$ in \eqref{eq:volume-velocity}. This is the key step to deal with the negative eigenvalues of $\Delta^L_g$, which produce oscillations.
\begin{lemma}\label{lem:volume-constraint}
In the coordinates \eqref{eq:slice-metric}, suppose
$\gamma=g_z+h\in\sliceball{\sigma}$, where $m\ge6$ and
$\sigma>0$ is sufficiently small. Then for fixed $c,C>0$,
\begin{equation}
 \label{eq:volume-taylor}
 \int_N R_\gamma\,dV_\gamma
 \ge c\|h_-\|_{L^2(N)}^2-C\|h_+\|_{\mathsf{H}^1(N)}^2,
\end{equation}
and,
\begin{equation}
 \label{eq:volume-velocity}
 |z'|+\|A\|_{L^2(N)}+\|h_-\|_{\mathsf{H}^1(N)}
 \le C\bigl(H+\|h_+\|_{\mathsf{H}^1(N)}\bigr).
\end{equation}
\end{lemma}

\begin{proof}
Let $E(g)\defeq\int_N R_g\,dV_g$. The first variation is
\begin{equation*}
 DE_g(k)=-\int_N
       \langle\Ric_g-\tfrac12R_gg,k\rangle_g\,dV_g.
\end{equation*}
At $g_0$, equation~\eqref{eq:deform} and $\delta_{g_0}h=0$ give
$DR_{g_0}(h)=\Delta\tr h$ and the second variation
\begin{equation}
 \label{eq:scalar-second-variation}
 D^2E_{g_0}(h,h)
 =-\tfrac12\langle\Delta^L_{g_0}h,h\rangle_{L^2}
     +\tfrac12\|d\tr h\|_{L^2}^2.
\end{equation}
Indeed, the Hessian term in $D\Ric$ has zero pairing with $h$ by integration by parts, while $\langle\Delta\tr h,\tr h\rangle=\|d\tr h\|^2$.
We expand $D^2 E_{g_0}(h,h)$ using $h=h_-+h_+$. If $W_-\ne\{0\}$,
choose an $L^2(N)$-orthonormal basis $e_1,\ldots,e_d$ of negative
eigentensors and write
\[
 \Delta^L_{g_0}e_\nu=\lambda_\nu e_\nu,\quad
 \lambda_\nu<0,\qquad
 h_-=\sum_{\nu=1}^d c_\nu e_\nu,\qquad
 \mu_- \defeq \min_\nu(-\lambda_\nu)>0 .
\]
Self-adjointness and spectral orthogonality imply
\[
 \langle\Delta^L_{g_0}h_-,h_+\rangle_{L^2}
 =\langle h_-,\Delta^L_{g_0}h_+\rangle_{L^2}=0.
\]
Since $\operatorname{tr}_{g_0}h_-=0$, we also have
$d\operatorname{tr}_{g_0}h=d\operatorname{tr}_{g_0}h_+$.
Thus \eqref{eq:scalar-second-variation} becomes
\begin{equation}
 \label{eq:scalar-spectral-expansion}
 D^2E_{g_0}(h,h)
 =\sum_{\nu=1}^d-\frac {\lambda_\nu}{2} c_\nu^2
   -\tfrac12\|\nabla h_+\|_{L^2(N)}^2+\int_N\langle\operatorname{Rm}_{g_0}h_+,h_+\rangle\,dV_{g_0}+\tfrac12\|d\operatorname{tr}_{g_0}h_+\|_{L^2(N)}^2 .
\end{equation}
Here we used
$\Delta^L_{g_0}=\nabla^*\nabla-2\operatorname{Rm}_{g_0}$.
The first term is at least $\mu_-\|h_-\|_{L^2(N)}^2/2$,
the last is nonnegative, and
\[
\left|\int_N\langle\operatorname{Rm}_{g_0}h_+,h_+\rangle\,dV_{g_0}\right|
 \le \|\operatorname{Rm}_{g_0}\|_{L^\infty(N)}\|h_+\|_{L^2(N)}^2 .
\]
It follows that
\begin{equation}
 \label{eq:scalar-spectral-lower-bound}
 D^2E_{g_0}(h,h)
 \ge \tfrac{\mu_-}{2}\|h_-\|_{L^2(N)}^2 -C_+\|h_+\|_{\mathsf{H}^1(N)}^2 .
\end{equation}
If $W_-=\{0\}$, omit the sums and constants involving $\mu_-$. The lower bound then contains only the $h_+$ term.

Note that $E(g_z)=DE_{g_z}=0$, so Taylor expansion at $g_z$ gives
\begin{equation}
 \label{eq:scalar-taylor-remainder}
 \left|E(g_z+h)-\tfrac12D^2E_{g_0}(h,h)\right|=\left|\int_0^1(1-t)
       \bigl(D^2E_{g_z+th}-D^2E_{g_0}\bigr)(h,h)\,dt\right|\le C\varepsilon\|h\|_{\mathsf{H}^1}^2
\end{equation}
when $\|g_z+th-g_0\|_{C^2}\le C\varepsilon$ for $0\le t\le1$.
Indeed, Lemma~\ref{lem:ricci-flat-decomposition} gives
$\|g_z-g_0\|_{C^2}\le C|z|
\le C\|\gamma-g_0\|_{C^2}\le C\sigma$.
Since $g_z+th=(1-t)g_z+t\gamma$, shrinking $\sigma$ ensures the estimate is uniform for all $0\le t\le 1$.

To verify the estimate, we first compute the second variation of $E$ at a general metric $g$ as follows:
 \begin{align*}
 D^2E_g(h,h)&=\int_N\Bigl\{
 -\tfrac12|\nabla h|_g^2+(\nabla_kh_{ij})(\nabla^ih^{kj})+\langle\delta_g h,d\operatorname{tr}_g h\rangle_g+\tfrac12|d\operatorname{tr}_g h|_g^2\\
 &\qquad+2\operatorname{Ric}_{g,ij}h^i{}_k h^{kj}-(\operatorname{tr}_gh)\langle\operatorname{Ric}_g,h\rangle_g-\tfrac12R_g|h|_g^2+\tfrac14R_g(\operatorname{tr}_g h)^2\Bigr\}\,dV_g\\
 & \defeq I_{1,g}+I_{2,g}.
 \end{align*}
Here $I_{1,g}$ and $I_{2,g}$ are the sums of the first and last four integrals, respectively. Each of the first four integrands contracts two copies of $\nabla^g h$ with three inverse metrics. We can estimate the integral difference w.r.t. $g$ and $g_0$ as
$$|I_{1,g}-I_{1,g_0}|\le C\|g-g_0\|_{C^1}\|h\|_{\mathsf{H}^1}^2.$$
The last four integrals vanish at $g_0$ and since $\|\operatorname{Ric}_g\|_{C^0}+\|R_g\|_{C^0} \le C\|g-g_0\|_{C^2}$, we have the estimate 
\[
|I_{2,g}-I_{2,g_0}|=|I_{2,g}|\le  C\|g-g_0\|_{C^2}\|h\|_{L^2}^2.
\]
Adding the estimates for $I_{1,g}$ and $I_{2,g}$ gives
\[
 |(D^2E_g-D^2E_{g_0})(h,h)|
 \le C\|g-g_0\|_{C^2}\|h\|_{\mathsf{H}^1}^2.
\]
Setting $g=g_z+th$, multiplying by $1-t$, and integrating over $[0,1]$ proves \eqref{eq:scalar-taylor-remainder}.

On the finite dimensional space
$W_-$ there is $C_->0$ with
$\|h_-\|_{\mathsf{H}^1(N)}^2\le C_-\|h_-\|_{L^2(N)}^2$.
Therefore the right hand side of \eqref{eq:scalar-taylor-remainder} satisfies
\[
 C\varepsilon\|h\|_{\mathsf{H}^1(N)}^2
 \le 2C\varepsilon\left(
       C_-\|h_-\|_{L^2(N)}^2+
       \|h_+\|_{\mathsf{H}^1(N)}^2\right).
\]
Choose $\varepsilon$ so that $2CC_-\varepsilon\le\mu_-/8$.
Combining this with \eqref{eq:scalar-spectral-lower-bound} gives
\[
 E(g_z+h)
 \ge \tfrac{\mu_-}{8}\|h_-\|_{L^2(N)}^2
       -(C_++2C\varepsilon)\|h_+\|_{\mathsf{H}^1(N)}^2 ,
\]
which proves \eqref{eq:volume-taylor}.

Integrating the scalar Gauss equation in \eqref{eq:volume-evolution} gives
\begin{equation*}
 \int_N |A|^2\,dV_\gamma+E(\gamma)=\Vol(N,\gamma)H^2.
\end{equation*}
Consequently, using $\gamma=g_z+h$, we obtain
\begin{equation}
 \label{eq:volume-second-form}
 \|A\|_{L^2(N)}^2+c\|h_-\|_{L^2(N)}^2
 \le C\bigl(H^2+\|h_+\|_{\mathsf{H}^1(N)}^2\bigr).
\end{equation}

To finish the proof of \eqref{eq:volume-velocity}, it remains to
estimate \(|z'|\). The quantitative shift estimate
\eqref{eq:slice-shift-low-order} in Lemma~\ref{lem:quantitative-shift}
and the uniform bound for \(u\) give
\[
 \|X\|_{\mathsf{H}^1(N)}\le C\|uA\|_{L^2(N)}
 \le C\|A\|_{L^2(N)}.
\]
Since Lemma~\ref{lem:ricci-flat-decomposition} and $P_{K_0} h=0$ give $P_{K_0}\gamma'=z'$,
projection of the first equation in \eqref{eq:volume-evolution} onto
$\ker\Delta^L_{g_0}$ gives
\begin{equation*}
 z'=P_{K_0}(2uA+\mathcal L_X\gamma).
\end{equation*}
The projection of $\mathcal L_Xg_0$ vanishes. The remaining Lie derivative satisfies 
$$\|\mathcal L_X (\gamma-g_0)\|_{L^2}\le C\|\gamma-g_0\|_{C^1}\|X\|_{\mathsf{H}^1},$$ proving
$|z'|\le C\|A\|_{L^2}$.
Together with \eqref{eq:volume-second-form} and finite dimensional norm equivalence on $W_-$, this proves \eqref{eq:volume-velocity}.
The same proof gives $|H|$ in place of $H$ without fixing orientation.
\end{proof}

\subsection{The transverse estimates}
The purpose of this section is to prove some estimates for $h_+$, as in previous sections we used $h_+$ to control other components in the metric tensor. Write $\dot g_zv=Dg_z[v]$ and $\ddot g_z(v,w)=D^2g_z[v,w]$. Recall that $P_{K_0}\dot g_z v=v$ and
$\Delta^L_{g_0}|_{W_+}$ has a positive spectral gap. If $W_+=\{0\}$, the estimate for $h_+$ below is immediate since $h_+=0$. We thus assume $W_+\neq \{0\}$.

Fix a length $L>0$. The following hypotheses and notation apply to the next three lemmas. Suppose a Ricci flat metric is written on $[a,a+(J+1)L]\times N$ as
\begin{equation*}
 g=u^2dr^2+\gamma_{ij}(dx^i+X^i dr)(dx^j+X^j dr),\qquad
 \gamma=g_z+h,
\end{equation*}
where the slices have constant mean curvature $H(r)$,
\[
 \fint_N u\,dV_\gamma=1,\qquad \delta_{g_0}(\gamma-g_0)=0,
\]
and $X(r)\perp\operatorname{Kill}(g_0)$ in $L^2(N)$.
Put $I_j=[a+jL,a+(j+1)L]$ for $0\le j\le J$.
Assume $P_{K_0}h=0$ and, for some $0<\varepsilon\le\varepsilon_0$,
\begin{equation}
 \label{eq:transverse-coefficient-smallness}
\max_{0\le j\le J}\bigl(
 \|\gamma-g_0\|_{C^4(I_j\times N)}
 +\|u-1\|_{C^4(I_j\times N)}+\|X\|_{C^4(I_j\times N)}\bigr)
 \le\varepsilon.
\end{equation}
 Define
\begin{equation*}
 e_j=\|h_+\|_{\mathsf{H}^2(I_j\times N)},\qquad
 H_j=\sup_{I_j}|H|,\qquad 0\leq j\leq J.
\end{equation*}
Here $\varepsilon_0>0$ is sufficiently small, depending only on a fixed neighborhood of $g_0$ and on $L$.

The first lemma gives estimates of the shift and normal speed in terms of the metric derivative, and more importantly, the estimates of the negative spectral term $h_-$ and the kernel term $z'$ in terms of mean curvature and the positive spectral term. Constants below are uniform on the fixed small neighborhood. Sobolev norms on the cross section at metrics uniformly close to $g_0$ are uniformly equivalent, and the fixed projections are bounded by the elliptic formula \eqref{eq:slice-projections} and the finite rank of $P_{K_0},P_-$.

\begin{lemma}\label{lem:transverse-preliminary}
Under the above hypotheses, on each $I_j$,
\begin{equation}
 \label{eq:transverse-velocity}
 |z'|+\|A\|_{L^2(N)}+\|h_-\|_{\mathsf{H}^1(N)}
 \leq C\bigl(|H|+\|h_+(r)\|_{\mathsf{H}^1(N)}\bigr)
 \leq C(H_j+e_j),
\end{equation}
where $A$ is the second fundamental form of leaves. The shift satisfies
\begin{equation}
 \label{eq:transverse-shift-h1}
 \|X\|_{\mathsf{H}^1(N)}
 \le C\varepsilon\|\gamma'\|_{L^2(N)}.
\end{equation}
\begin{equation}
 \label{eq:transverse-shift-h2}
 \|X\|_{\mathsf{H}^2(N)}
 \le C\varepsilon\|\gamma'\|_{\mathsf{H}^1(N)}.
\end{equation}
\begin{equation}
 \label{eq:transverse-shift-r}
 \|X'\|_{\mathsf{H}^1(N)}
 \le C\varepsilon\bigl(
   \|\gamma''\|_{L^2(N)}+\|\gamma'\|_{\mathsf{H}^1(N)}\bigr).
\end{equation}
The normal speed satisfies
\begin{equation}
 \label{eq:transverse-normal-speed-bound}
 \|u-1\|_{\mathsf{H}^2(N)}
 \le C\left\|u|A|^2- \fint_N u|A|^2\,dV_\gamma\right\|_{L^2(N)}\le C\varepsilon\|A\|_{L^2(N)}
 \le C\varepsilon\|\gamma'\|_{L^2(N)}.
\end{equation}
For every fixed integer $s$,
\begin{equation}
 \label{eq:transverse-negative-first}
 \|h_-\|_{L^2(I_j;\mathsf{H}^s(N))}
 +\|h_-'\|_{L^2(I_j;\mathsf{H}^s(N))}
 +\|z'\|_{L^2(I_j)}
 \le C_s(H_j+e_j).
\end{equation}
\end{lemma}

\begin{proof}
Lemma~\ref{lem:volume-constraint}, with $|H|$ in place of $H$, gives the first inequality in
\eqref{eq:transverse-velocity}. The uniform trace estimate
\begin{equation}\label{eq:trace-estimate_plus}
 \sup_{r\in I_j}\|h_+(r)\|_{\mathsf{H}^1(N)}
 \le C_L\|h_+\|_{\mathsf{H}^2(I_j\times N)}
\end{equation}
gives the second.

We next estimate the shift and normal speed.

The contracted Codazzi identity $\delta_\gamma A=0$ and $\gamma'=2uA+\mathcal L_X\gamma$ in \eqref{eq:volume-evolution} give
\[
 \delta_\gamma\bigl(u^{-1}(\gamma'-\mathcal L_X\gamma)\bigr)=0.
\]
Differentiating the fixed slice condition
$\delta_{g_0}(\gamma-g_0)=0$ gives $\delta_{g_0}\gamma'=0$. Recall $\mathcal B_{g_0}$ from \eqref{eq:fixed-slice-operator}. The above two identities give the decomposition 
\begin{align*}
 \mathcal B_{g_0}X
 =\bigl(\mathcal B_{g_0}X-\delta_\gamma(u^{-1}\mathcal L_X\gamma)\bigr)+\bigl(\delta_\gamma(u^{-1}\gamma')-\delta_{g_0}\gamma'\bigr).
\end{align*}

 For
$Y\perp\operatorname{Kill}(g_0)$,
Remark~\ref{rem:fixed-slice-low-order} gives
\begin{equation}\label{eq:H-minusone-to-H-one}
 \|Y\|_{\mathsf{H}^1}
 \le C\|\mathcal B_{g_0}Y\|_{\mathsf{H}^{-1}},
\end{equation}
and the elliptic estimate for the fixed operator also gives
\begin{equation}\label{eq:L-two-to-H-two}
 \|Y\|_{\mathsf{H}^2}
 \le C\bigl(\|\mathcal B_{g_0}Y\|_{L^2}+\|Y\|_{L^2}\bigr)
 \le C\|\mathcal B_{g_0}Y\|_{L^2}.
\end{equation}

We estimate $\mathcal B_{g_0}X$ by dealing with the two differences in the decomposition. The above discussion then gives the estimates of $X$. We start with the easier second difference. When the divergences are expressed using covariant derivatives with respect to $g_0$, the coefficients of the first order terms in $\delta_\gamma(u^{-1}\gamma')-\delta_{g_0}\gamma'$ involve $u^{-1}\gamma^{ij}-g_0^{ij}$. The remaining coefficients involve first derivatives of $u$ w.r.t. variables on the cross section and $\gamma-g_0$. By \eqref{eq:transverse-coefficient-smallness}, the former are $O(\varepsilon)$ in $C^1$ and the latter are $O(\varepsilon)$ in $C^0$. Integration by parts against a test one-form $\omega$ gives
\begin{align*}
 \bigl|\langle\delta_\gamma(u^{-1}\gamma')-\delta_{g_0}\gamma',\omega
 \rangle\bigr|\le C\bigl(\|\gamma-g_0\|_{C^1}+\|u-1\|_{C^1}\bigr)\|\gamma'\|_{L^2}\|\omega\|_{\mathsf{H}^1}
 \le C\varepsilon\|\gamma'\|_{L^2}\|\omega\|_{\mathsf{H}^1}.
\end{align*}
Hence
\[
 \|\delta_\gamma(u^{-1}\gamma')-\delta_{g_0}\gamma'\|_{\mathsf{H}^{-1}}
 \le C\varepsilon\|\gamma'\|_{L^2}.
\]
For the first difference, use
\begin{align*}
 \mathcal B_{g_0}X-\delta_\gamma(u^{-1}\mathcal L_X\gamma)
 ={}&\delta_{g_0}(\mathcal L_X\gamma)
      -\delta_\gamma(u^{-1}\mathcal L_X\gamma)-\delta_{g_0}\mathcal L_X(\gamma-g_0).
\end{align*}
The same integration by parts applies with
$\mathcal L_X\gamma$ in place of $\gamma'$. Moreover, 
\begin{align*}
 \bigl(\mathcal L_X(\gamma-g_0)\bigr)_{ij}
 ={}&X^k\nabla^{g_0}_k(\gamma-g_0)_{ij}+(\gamma-g_0)_{kj}\nabla^{g_0}_iX^k
   +(\gamma-g_0)_{ik}\nabla^{g_0}_jX^k,
\end{align*}
so
\[
\|\mathcal L_X(\gamma-g_0)\|_{L^2}
\le C\|\gamma- g_0\|_{C^1}\|X\|_{\mathsf{H}^1}
\]

It follows that
\begin{align*}
 \|\mathcal B_{g_0}X-\delta_\gamma(u^{-1}\mathcal L_X\gamma)\|_{\mathsf{H}^{-1}}
 \le
 C\varepsilon\|\mathcal L_X\gamma\|_{L^2}+C\|\mathcal L_X(\gamma-g_0)\|_{L^2}
 \le C\varepsilon\|X\|_{\mathsf{H}^1}.
\end{align*}
Combining the two difference estimates gives the bound on $ \|\mathcal B_{g_0}X\|_{\mathsf{H}^{-1}}$.

Applying \eqref{eq:H-minusone-to-H-one} therefore gives
\[
 \|X\|_{\mathsf{H}^1}
 \le C\|\mathcal B_{g_0}X\|_{\mathsf{H}^{-1}}\le
C\varepsilon\|X\|_{\mathsf{H}^1}+C\varepsilon\|\gamma'\|_{L^2}.
\]
After decreasing $\varepsilon_0$ so that $C\varepsilon_0\le1/2$, absorption yields \eqref{eq:transverse-shift-h1}.
To obtain the $\mathsf{H}^2$ bound of $X$, estimate the same divergence difference in $L^2$, 
\[
 \|\delta_\gamma(u^{-1}\gamma')-\delta_{g_0}\gamma'\|_{L^2}
 \le C\varepsilon\|\gamma'\|_{\mathsf{H}^1},
\] 
and use the same Lie derivative bound in $\mathsf{H}^1$:
\[
\|\mathcal L_X(\gamma-g_0)\|_{\mathsf{H}^1}
 \le C\|\gamma-g_0\|_{C^2}\|X\|_{\mathsf{H}^2}.
\]
Consequently, the same decomposition of $\mathcal B_{g_0}X$ gives
\begin{align*}
 \|\mathcal B_{g_0}X
       -\delta_\gamma(u^{-1}\mathcal L_X\gamma)\|_{L^2}
 &\le C\varepsilon\|\mathcal L_X\gamma\|_{\mathsf{H}^1} +C\|\mathcal L_X(\gamma-g_0)\|_{\mathsf{H}^1}\le C\varepsilon\|X\|_{\mathsf{H}^2}.
\end{align*}
The estimate \eqref{eq:L-two-to-H-two} now gives
\[
 \|X\|_{\mathsf{H}^2}
 \le  C\varepsilon\|X\|_{\mathsf{H}^2}+C\varepsilon\|\gamma'\|_{\mathsf{H}^1}.
\]
Another absorption proves \eqref{eq:transverse-shift-h2}.

For the radial derivative, differentiate the contracted Codazzi
identity to obtain
\begin{align*}
 \delta_\gamma(u^{-1}\mathcal L_{X'}\gamma)
 ={}&\delta_\gamma(u^{-1}\gamma'')
   +(\partial_r\delta_\gamma)
       \bigl(u^{-1}(\gamma'-\mathcal L_X\gamma)\bigr)-\delta_\gamma
       \bigl(u^{-2}u'(\gamma'-\mathcal L_X\gamma)\bigr)
   -\delta_\gamma(u^{-1}\mathcal L_X\gamma').
\end{align*}
Here $\partial_r\delta_\gamma$ differentiates only the coefficients
of the divergence operator. The mixed derivative bounds in
\eqref{eq:transverse-coefficient-smallness} and the same integration
by parts give
\begin{align*}
 \bigl\|&(\partial_r\delta_\gamma)
       \bigl(u^{-1}(\gamma'-\mathcal L_X\gamma)\bigr)
       \bigr\|_{\mathsf{H}^{-1}}+\bigl\|\delta_\gamma
       \bigl(u^{-2}u'(\gamma'-\mathcal L_X\gamma)\bigr)
       \bigr\|_{\mathsf{H}^{-1}}\\&\le C\bigl(\|\gamma'\|_{C^1}+\|u'\|_{C^0}\bigr)\bigl(\|\gamma'\|_{L^2}+\|X\|_{\mathsf{H}^1}\bigr)\\&\le C\varepsilon
       \bigl(\|\gamma'\|_{L^2}+\|X\|_{\mathsf{H}^1}\bigr).
\end{align*}
The remaining term satisfies
\[
 \|\delta_\gamma(u^{-1}\mathcal L_X\gamma')\|_{\mathsf{H}^{-1}}
 \le C\|\gamma'\|_{C^1}\|X\|_{\mathsf{H}^1}
 \le C\varepsilon\|X\|_{\mathsf{H}^1}.
\]
Differentiating the fixed slice condition twice gives
$\delta_{g_0}\gamma''=0$. Applying the preceding divergence and
Lie derivative estimates with $X'$ in place of $X$ and $\gamma''$
in place of $\gamma'$ therefore gives
\[
 \|\mathcal B_{g_0}X'\|_{\mathsf{H}^{-1}}
 \le C\varepsilon\bigl(
   \|X'\|_{\mathsf{H}^1}+\|\gamma''\|_{L^2}
   +\|\gamma'\|_{L^2}+\|X\|_{\mathsf{H}^1}\bigr).
\]
The fixed normalization implies
$X'\perp\operatorname{Kill}(g_0)$, so the same inverse estimate \eqref{eq:H-minusone-to-H-one} applies to $X'$. Absorbing its $\mathsf{H}^1$ norm on the right,
using \eqref{eq:transverse-shift-h1}, and bounding
$\|\gamma'\|_{L^2}$ by $\|\gamma'\|_{\mathsf{H}^1}$ proves \eqref{eq:transverse-shift-r}.
Here $\|\gamma'\|_{\mathsf{H}^1(N)}$ controls $\gamma'$ and
$\nabla_{g_0}\gamma'$ on a leaf, and does not include the radial
derivative $\gamma''$.

The mean curvature variation equation and its integration over $N$ are
\begin{equation}
 \label{eq:transverse-normal-speed}
 H'=\Delta_\gamma u-u|A|^2,\qquad
 H'=-\fint_N u|A|^2\,dV_\gamma.
\end{equation}
The second equality uses that $H'=H'(r)$ is independent of cross section variables and the integral of a Laplacian on a closed manifold vanishes. Subtracting the second equality from the first in
\eqref{eq:transverse-normal-speed} gives the mean-zero
elliptic equation
\begin{equation*}
 \Delta_\gamma(u-1)
 =u|A|^2-\fint_N u|A|^2\,dV_\gamma,\qquad
 \fint_N (u-1)\,dV_\gamma=0.
\end{equation*}
The coefficient bounds \eqref{eq:transverse-coefficient-smallness}
give uniform ellipticity, bounded coefficients, and a positive lower
bound for the first nonzero scalar eigenvalue. The mean zero
elliptic estimate and the smallness of $A$ therefore give \eqref{eq:transverse-normal-speed-bound}.
Here the last inequality follows directly from $2uA=\gamma'-\mathcal L_X\gamma$ and
\eqref{eq:transverse-shift-h1} as follows. First
$$\|A\|_{L^2}\le C(\|\gamma'\|_{L^2}+\|X\|_{\mathsf{H}^1})
\le C\|\gamma'\|_{L^2}.$$ 
Conversely,
$$\|\gamma'\|_{L^2}\le C\|A\|_{L^2}+C\|X\|_{\mathsf{H}^1}
\le C\|A\|_{L^2}+C\varepsilon\|\gamma'\|_{L^2},$$ so absorption gives
$$\|\gamma'\|_{L^2}\le C\|A\|_{L^2}.$$

Finally,
$h_-'=P_-\gamma'-P_-\dot g_z z'$.
The bound $\|\gamma'\|_{L^2}\le C\|A\|_{L^2}$ above and finite
dimensionality of $W_-$ therefore imply \eqref{eq:transverse-negative-first} for every fixed integer $s$.
Here, the $L^2$ norm is taken on the interval $I_j$ and the $\mathsf{H}^s$ norm is taken on the cross section $N$.
\end{proof}

The second lemma is an intermediate technical step to study the equation satisfied by the positive spectral term $h_+$.
\begin{lemma}\label{lem:transverse-equation}
Under the above hypotheses, there is a tensor $\mathcal E$ such that, on every subinterval $I$,
\begin{equation}
 \label{eq:transverse-projected-equation}
 \begin{gathered}
 P_S(-\gamma''+\Delta^L_{g_0}h)=\mathcal E,\\
 \|\mathcal E\|_{L^2(I\times N)}\le C\varepsilon \bigl(\|h\|_{\mathsf{H}^2(I\times N)}+\|z'\|_{L^2(I)}+\|z''\|_{L^2(I)}\bigr).
 \end{gathered}
\end{equation}
On each $I_j$,
\begin{equation}
 \label{eq:transverse-negative-control}
 \begin{aligned}
 \|h_-\|_{\mathsf{H}^2(I_j\times N)}&\le C(H_j+e_j),\\
 \|z''\|_{L^2(I_j)}+\|\mathcal E\|_{L^2(I_j\times N)}
      &\le C\varepsilon(H_j+e_j).
 \end{aligned}
\end{equation}
The tensor $h_+$ satisfies
\begin{equation}
 \label{eq:transverse-elliptic-equation}
 \begin{gathered}
 (-\partial_r^2+\Delta^L_{g_0})h_+=\mathcal F_+,\qquad
 \mathcal F_+\defeq P_+\mathcal E+P_+\dot g_z z''+P_+\ddot g_z(z',z'),\\
 \|\mathcal F_+\|_{L^2(I_j\times N)}
       \le C\varepsilon(e_j+H_j).
 \end{gathered}
\end{equation}
\end{lemma}

\begin{proof}
We first obtain an elliptic equation for $h$.
By \eqref{eq:deform}, at $g_0$,
\begin{equation*}
 2P_S D\operatorname{Ric}_{g_0}(h)=\Delta^L_{g_0}h,
 \quad h\in S,
\end{equation*}
because the remaining Hessian is orthogonal to $S$. $\operatorname{Ric}(g_z)=0$ and Taylor expansion at $g_z$ imply
\begin{equation}
 \label{eq:transverse-ricci-expansion}
 \|2P_S\operatorname{Ric}(g_z+h)-\Delta^L_{g_0}h\|_{L^2(N)}
 \leq C\varepsilon\|h\|_{\mathsf{H}^2(N)}.
\end{equation}

Write $D_r\defeq \partial_r-\mathcal L_X$. Since
$\partial_r-X=u\nu$, the normal variation formula proved in
Proposition~\ref{prop:cmc-volume} gives
\begin{equation*}
 D_r\gamma=2uA,\qquad
 D_rA=-\nabla_\gamma^2u
       +u\bigl(\operatorname{Ric}(\gamma)-HA+2A\circ A\bigr).
\end{equation*}
Moreover,
\begin{equation*}
 D_r^2\gamma
 =\gamma''-\mathcal L_{X'}\gamma-2\mathcal L_X\gamma'
      +\mathcal L_X^2\gamma.
\end{equation*}
Applying $D_r$ to $D_r\gamma=2uA$ and using these identities yields
\begin{align*}
 \gamma''={}&2u^2\operatorname{Ric}(\gamma)-2u\nabla_\gamma^2u
       +4u^2A\circ A-2u^2HA+\mathcal L_{X'}\gamma+2\mathcal L_X\gamma'
       -\mathcal L_X^2\gamma+2(u'-X(u))A.
\end{align*}
Set
\begin{equation*}
 \mathcal R=2P_S\operatorname{Ric}(g_z+h)-\Delta^L_{g_0}h.
\end{equation*}
By \eqref{eq:transverse-ricci-expansion},
$\|\mathcal R\|_{L^2}\le C\varepsilon\|h\|_{\mathsf{H}^2}$. Applying $P_S$
to the formula for $\gamma''$ and moving its right-hand side to the left
defines
\begin{align*}
 \mathcal E={}&-\mathcal R-2P_S\bigl((u^2-1)\operatorname{Ric}(\gamma)\bigr)
       +2P_S(u\nabla_\gamma^2u)-4P_S(u^2A\circ A)\\
 &+2P_S(u^2HA)-P_S\mathcal L_{X'}\gamma
       -2P_S\mathcal L_X\gamma'+P_S\mathcal L_X^2\gamma
       -2P_S\bigl((u'-X(u))A\bigr).
\end{align*}
Thus $P_S(-\gamma''+\Delta^L_{g_0}h)=\mathcal E$. The identities
$P_S\nabla_{g_0}^2f=0$ and $P_S\mathcal L_Yg_0=0$ give
\begin{equation*}
 \|P_S(u\nabla_\gamma^2u)\|_{L^2}
       \le C\varepsilon\|u-1\|_{\mathsf{H}^2},\qquad
 \|P_S\mathcal L_{X'}\gamma\|_{L^2}
       \le C\varepsilon\|X'\|_{\mathsf{H}^1}.
\end{equation*}
Also $\|\operatorname{Ric}(g_z+h)\|_{L^2}\le C\|h\|_{\mathsf{H}^2}$,
and each remaining term contains a small coefficient controlled by
the assumption \eqref{eq:transverse-coefficient-smallness}. Using
\eqref{eq:transverse-shift-h1}, \eqref{eq:transverse-shift-h2},
\eqref{eq:transverse-shift-r}, and
\eqref{eq:transverse-normal-speed-bound} from Lemma~\ref{lem:transverse-preliminary} therefore gives \eqref{eq:transverse-projected-equation} on every subinterval $I$.
Indeed, $\gamma'=\dot g_z z'+h'$ and
$\gamma''=\dot g_z z''+\ddot g_z(z',z')+h''$. The $C^4(I_j\times N)$ bound \eqref{eq:transverse-coefficient-smallness} gives $|z'|\le C\varepsilon$, so the quadratic term is bounded by
$C\varepsilon|z'|$.

We next control the second radial derivative of $h_-$. Since $P_{K_0}(g_z-g_0)=z$ and $P_{K_0} h=0$, differentiation gives $P_{K_0}\dot g_z=\operatorname{Id}$ and $P_{K_0}\ddot g_z=0$. Projecting \eqref{eq:transverse-projected-equation} onto $K_0$ and $W_-$ gives
\begin{equation}\label{eq:projected-equations}
 z''=-P_{K_0}\mathcal E,
 \quad
 h_-''=\Delta^L_{g_0}h_--P_-\mathcal E-P_-\dot g_z z''-P_-\ddot g_z(z',z').
\end{equation}
Here $P_-\dot g_z=O(|z|)=O(\varepsilon)$, $\ddot g_z$ is bounded,
and $|z'|\le C\varepsilon$. Finite dimensional norm equivalence and
\eqref{eq:transverse-negative-first} imply
\begin{equation}\label{eq:H^2-bound-of-h-minus}
 \|h_-\|_{\mathsf{H}^2(I_j\times N)}
       \le C(H_j+e_j)+C\|h_-''\|_{L^2(I_j\times N)},
\end{equation}
and the two projected equations \eqref{eq:projected-equations} imply
\begin{equation}\label{eq:second-derivative-bound}
 \|z''\|_{L^2(I_j)}+\|h_-''\|_{L^2(I_j\times N)}
       \le C(H_j+e_j)+C\|\mathcal E\|_{L^2(I_j\times N)}.
\end{equation}
Since $\|h\|_{\mathsf{H}^2}\le e_j+\|h_-\|_{\mathsf{H}^2}$ and
$\|z'\|_{L^2}\le C(H_j+e_j)$, substitution into \eqref{eq:transverse-projected-equation} gives
\begin{equation*}
 \|\mathcal E\|_{L^2(I_j\times N)}
 \le C\varepsilon\bigl(H_j+e_j+\|\mathcal E\|_{L^2(I_j\times N)}\bigr).
\end{equation*}
After decreasing $\varepsilon_0$ so that $C\varepsilon_0<1/2$, the
last term is absorbed. Substitution into the bounds for \eqref{eq:H^2-bound-of-h-minus} and \eqref{eq:second-derivative-bound} proves \eqref{eq:transverse-negative-control}.

Finally apply $P_+$ to \eqref{eq:transverse-projected-equation}.
Since $P_+\dot g_z=O(\varepsilon)$, we obtain \eqref{eq:transverse-elliptic-equation}.
The last bound uses \eqref{eq:transverse-negative-control},
\eqref{eq:transverse-negative-first}, and $|z'|=O(\varepsilon)$.
\end{proof}

The third lemma gives a kind of convexity of the $L^2$ norm of $h_+$, which can be viewed as a variant of the three circle type estimates. This is the key to exponential decay.
\begin{lemma}\label{lem:transverse-estimate}
Under the above hypotheses, there are $C<\infty$ and $0<\theta<1$, depending only on the fixed neighborhood of $g_0$ and on $L$, such that
\begin{equation}
 \label{eq:transverse-green}
 e_j\leq C\bigl(\theta^je_0+\theta^{J-j}e_J\bigr)
       +C\varepsilon\sum_{k=0}^J\theta^{|j-k|}(H_k+e_k).
\end{equation}
The constants do not depend on $a$, $J$, or the length of the coordinate interval.
\end{lemma}

\begin{proof}
We apply the equation from Lemma~\ref{lem:transverse-equation}. Put
\[
 \mathscr P=-\partial_r^2+\Delta^L_{g_0}|_{W_+}.
\]
Let $\{e_\nu\}$ be an orthonormal basis of $W_+$ with
$\Delta^L_{g_0}e_\nu=\lambda_\nu e_\nu$ and
$\lambda_\nu\ge\lambda_*>0$. For $f\in L^2(\mathbb R;W_+)$,
write $f(r)=\sum_\nu f_\nu(r)e_\nu$ and Fourier transform in $r$.
The solution $v=\mathscr P^{-1}f$ satisfies
\[
 \widehat v_\nu(\tau)
 =\frac{\widehat f_\nu(\tau)}{\tau^2+\lambda_\nu},\qquad
 \sum_\nu\int_{\mathbb R}
 (1+\tau^2+\lambda_\nu)^2|\widehat v_\nu(\tau)|^2\,d\tau
 \le C\|f\|_{L^2(\mathbb R\times N)}^2 .
\]
The expression on the left is equivalent to
$\|v\|_{\mathsf{H}^2(\mathbb R\times N)}^2$ by elliptic regularity
on $N$. Thus
\begin{equation}\label{eq:global-inequality}
\|\mathscr P^{-1}f\|_{\mathsf{H}^2(\R\times N)}\le C\|f\|_{L^2(\R\times N)}.
\end{equation}
 
The inverse has the spectral kernel representation
\[
 (\mathscr P^{-1}f)(r)
 =\int_{\mathbb R}K(r,s)f(s)\,ds,\qquad
 K(r,s) \defeq \frac12{e^{-|r-s|({\Delta^L_{g_0}|_{W_+}})^{\frac12}}}{({\Delta^L_{g_0}|_{W_+}})^{-\frac12}}.
\]

Suppose $f$ is supported in $I_k$. If $|j-k|\ge2$, then $|r-s|\ge(|j-k|-1)L$ for $r\in I_j$ and $s\in I_k$.
For $d=0,1,2$, the positive spectral gap gives
\[
 \|\partial_r^d K(r,s)\|_
 {L^2(N)\to\mathsf{H}^{2-d}(N)}
 \le C_L e^{-c|r-s|},\quad|r-s|\ge L,
\]
where $\mathsf{H}^0=L^2$. Indeed, up to constants depending on $\lambda_*$, the spectral multiplier on the left is bounded by $\lambda^{1/2}e^{-|r-s|\sqrt\lambda}$. Maximizing over $\lambda\ge\lambda_*$ gives $$\lambda^{1/2}e^{-|r-s|\sqrt\lambda}\le C_Le^{-c|r-s|}, \quad |r-s|\ge L.$$
Integrating first in $s$ and then in $r$ over intervals of length $L$ proves the estimate \eqref{eq:transverse-resolvent} for $|j-k|\ge2$. When $|j-k|\le1$, use the global estimate \eqref{eq:global-inequality} and enlarge its constant. In conclusion, for some $\theta\in(0,1)$ depending only on $L$ and $\lambda_*$, we obtain
\begin{equation}
 \label{eq:transverse-resolvent}
 \|\mathscr P^{-1}f\|_{\mathsf{H}^2(I_j\times N)}
 \le C\theta^{|j-k|}\|f\|_{L^2(I_k\times N)} .
\end{equation}

For $J\ge2$, set $b=a+(J+1)L$ and choose a smooth cutoff
$\chi:[a,b]\to[0,1]$ which vanishes near $a,b$, equals $1$ on $[a+L,b-L]$, and satisfies
$\|\chi'\|_\infty+\|\chi''\|_\infty\le C_L$.
Extend $v\defeq \chi h_+$ by zero to $\mathbb R\times N$. Equation
\eqref{eq:transverse-elliptic-equation} gives
\[
 \mathscr Pv=\chi \mathcal F_+ + c_0+c_J,\quad
 c_0+c_J=-\chi''h_+-2\chi'\partial_rh_+,
\]
where $c_0$ and $c_J$ are supported in $I_0\times N$ and $I_J\times N$, respectively. The definitions of $e_0$ and $e_J$ give
\[
 \|c_0\|_{L^2(I_0\times N)}\le C_Le_0,\qquad
 \|c_J\|_{L^2(I_J\times N)}\le C_Le_J .
\]
Split $\chi \mathcal F_+$ into its restrictions to $I_k\times N$ and apply
\eqref{eq:transverse-resolvent} to each term. Since $v=h_+$ for $1\le j\le J-1$, this gives
\begin{align*}
 e_j&\le C\left(\theta^je_0+\theta^{J-j}e_J\right)+C\sum_{k=0}^J\theta^{|j-k|}
 \|\mathcal F_+\|_{L^2(I_k\times N)}\\
&\le C\left(\theta^je_0+\theta^{J-j}e_J\right)+C\varepsilon\sum_{k=0}^J\theta^{|j-k|}(H_k+e_k).
\end{align*}
The second inequality is
\eqref{eq:transverse-elliptic-equation}. For $j=0,J$ the same
inequality follows by choosing $C\ge1$, because its right side
already contains $Ce_j$. This also handles $J=0,1$.
We have proved \eqref{eq:transverse-green}. All cutoff widths, spectral gaps,
and elliptic constants are fixed independently of $a$ and $J$.
\end{proof}

\subsection{Continuation and decay}
By Proposition~\ref{prop:cmc-volume}, the function $\log\vol(N,\gamma(r))$ satisfies
\begin{equation}
 \label{eq:volume-concavity}
 \log\vol(N,\gamma(r))'=H(r)\geq0,\qquad
 H'=-\fint_N u|A|^2\,dV_\gamma,\qquad H\longrightarrow0.
\end{equation}
 We have
\begin{equation}
 \label{eq:volume-tail}
 \int_a^\infty H(r)\,dr=\log V_\infty-\log\vol(N,\gamma(a))\longrightarrow0.
\end{equation}

We recall the consequence of Lemma~\ref{lem:transverse-estimate}
that will be used. Fix the length $L$, put
$I_j=[a+jL,a+(j+1)L]$, and set
\begin{equation}
 \label{eq:volume-blocks}
 e_j=\|h_+\|_{\mathsf{H}^2(I_j\times N)},\qquad
 H_j=H(a+jL)=\sup_{I_j}|H|,\qquad Q_j=H_j+e_j.
\end{equation}
Lemma~\ref{lem:transverse-preliminary} bounds $|z'|$, $\|A\|_{L^2}$, and $\|h_-\|_{\mathsf{H}^1}$ by $CQ_j$.
Apply Lemma~\ref{lem:transverse-estimate} to any finite sequence of $I_j\times N$'s on which \eqref{eq:transverse-coefficient-smallness} holds.
By \eqref{eq:volume-concavity}, $H\ge0$ and $H'\le0$, so $\sup_{I_k}|H|=H(a+kL)=H_k$. Thus \eqref{eq:transverse-green} gives
\begin{equation}
 \label{eq:volume-transverse}
 e_j\le C\bigl(e_0\theta^j+e_J\theta^{J-j}\bigr)
       +C\varepsilon\sum_{k=0}^J\theta^{|j-k|}Q_k,
       \qquad 0\le j\le J.
\end{equation}
Fix $C,\theta$ on an initial slice neighborhood, enlarging $C$ to cover the estimates below. They remain valid on smaller neighborhoods.
The number $\varepsilon$ is as small as required by shrinking the neighborhood and then choosing the initial parameter $a$ sufficiently large. Estimate~\eqref{eq:slice-coefficient-smallness} gives the $C^4(I_j\times N)$ bound \eqref{eq:transverse-coefficient-smallness} required by Lemma~\ref{lem:transverse-estimate}.

Now we are in a position to prove the main theorem.
\begin{proof}[Proof of Theorem~\ref{thm:main}]
The first step is to show that the slice coordinates in Proposition \ref{prop:slice-path} can be extended to infinity. Fix $m\ge6$ and the balls $\sliceball{2\sigma}$ and
$\sliceball{\sigma}$. We construct coordinates in
$\sliceball{2\sigma}$ and carry out the estimates while
$\gamma(r)\in\sliceball{\sigma}$. Choose $\varepsilon>0$ so small that
\begin{equation*}
 C\varepsilon\frac{1+\theta}{1-\theta}<\frac14,
 \qquad \theta(1+2C\varepsilon)<1.
\end{equation*}
Shrink $\sigma$ as required for the fixed slice neighborhood,
and choose $\rho>0$ so that
\begin{equation*}
 |z|\le2\rho\quad\Longrightarrow\quad
 \|g_z-g_0\|_{C^{m,\alpha}(N)}<\sigma/4.
\end{equation*}
Shrink $\rho$ further if necessary and choose $a_0$ sufficiently
large that \eqref{eq:slice-coefficient-smallness} gives
\eqref{eq:transverse-coefficient-smallness} on every $I_j\times N$ with $\gamma(r)\in\overline{\sliceball{\sigma}}$
and $z(r)\in\overline{\kernelball{\rho}}$
and initial parameter $a\ge a_0$.
Proposition~\ref{prop:slice-path}, applied to
$\overline{\sliceball{\sigma}}\subset\sliceball{2\sigma}$,
has the following uniform consequence.
There is $\omega(a)\downarrow0$ as $a\to\infty$ such that, on
every normalized slice lift with metric representatives in $\overline{\sliceball{\sigma}}$ on an interval $I_*\subset[a,\infty)$,
\begin{equation}
 \label{eq:continuation-smallness}
 \begin{gathered}
 \|h(r)\|_{C^{m,\alpha}(N)}+\|A(r)\|_{L^2}
       \le\omega(a)\qquad(r\in I_*),\\
 \|h\|_{\mathsf{H}^2(I\times N)}\le\omega(a)
       \qquad(I\subset I_*,\ |I|=L),\\
 H(a)+\log V_\infty-\log\vol(N,\gamma(a))\le\omega(a).
 \end{gathered}
\end{equation}
Use the CMC construction and the local slice map to choose a leaf
along the prescribed asymptotic sequence with $a\ge a_0$ and
\begin{equation*}
 \|\gamma(a)-g_0\|_{C^{m,\alpha}(N)}<\sigma/4,\qquad
 |z(a)|<\rho/4 .
\end{equation*}
Let $T$ be the infimum for which
\begin{equation*}
 |z(T)|=\rho
 \quad\text{or}\quad
 \|\gamma(T)-g_0\|_{C^{m,\alpha}(N)}=\sigma.
\end{equation*}
By Proposition~\ref{prop:slice-path}, a normalized lift cannot terminate at a finite $T$ while $\gamma(r)$ remains in
$\overline{\sliceball{\sigma}}\subset\sliceball{2\sigma}$. $T$ is finite only when
$z(T)\in\partial\kernelball{\rho}$ or
$\gamma(T)\in\partial\sliceball{\sigma}$. If neither occurs, we can set $T=\infty$ which is the desired result. On $[a,T)$, summing \eqref{eq:volume-transverse} over any $J+1$ complete $I_j\times N$ and absorbing the sum containing $e_j$ gives
\begin{equation}
 \label{eq:volume-sum}
 \sum_{j=0}^J e_j
 \le C(e_0+e_J)+C\varepsilon\sum_{j=0}^JH_j.
\end{equation}
Since $H$ is nonnegative and decreasing,
\begin{equation*}
 \sum_{j=0}^JH_j
 \le H(a)+L^{-1}\int_a^{a+JL}H
 \le H(a)+L^{-1}(\log V_\infty-\log\vol(N,\gamma(a))).
\end{equation*}
Integrating \eqref{eq:volume-velocity} for each $I_j\times N$ through $j=0$ to $j=J$ yields
\begin{equation}
 \label{eq:volume-displacement}
 \int_a^{a+(J+1)L}|z'|
 \le C\bigl(e_0+e_J+H(a)+\log V_\infty-\log\vol(N,\gamma(a))\bigr)
 \le C\omega(a).
\end{equation}
For the terminal interval of length possibly less than $L$, use instead $|z'|\le C\|A\|_{L^2}$, proved in Lemma~\ref{lem:volume-constraint}. Its contribution is at most $CL\omega(a)$. Consequently, for every $t<T$,
\begin{equation*}
 |z(t)-z(a)|\le\int_a^t|z'|\le C_L\omega(a).
\end{equation*}
Choose $a\ge a_0$ sufficiently large that $C_L\omega(a)<\rho/4$ and $\omega(a)<\sigma/4$. Then throughout $[a,T)$,
\begin{equation*}
 \begin{gathered}
 |z(t)|<\rho/2,\\
 \|\gamma(t)-g_0\|_{C^{m,\alpha}(N)}
 \le
 \|g_{z(t)}-g_0\|_{C^{m,\alpha}(N)}+\|h(t)\|_{C^{m,\alpha}(N)}<\sigma/2.
 \end{gathered}
\end{equation*}
These strict bounds exclude $\partial\kernelball{\rho}$ and $\partial\sliceball{\sigma}$. If $T$ were finite, the recentered extension in
Proposition~\ref{prop:slice-path} would extend the same normalized lift past $T$, still with metric representatives in $\sliceball{\sigma}$, a contradiction.
Therefore $T=\infty$, establishing the first step.

Now that $T=\infty$, this enables us to take $J$ to infinity. The next step is to prove the decay of $H_j+e_j$ using a difference inequality. Because the global lift remains in $\overline{\sliceball{\sigma}}$, the convergence conclusion of Proposition~\ref{prop:slice-path}, applied at the centres $a+JL\to\infty$, gives $e_J\to0$. Moreover, monotonicity of $H$ and \eqref{eq:volume-tail} give
\begin{equation*}
 \sum_{j=0}^\infty H_j
 \le H(a)+L^{-1}\int_a^\infty H(r)\,dr<\infty.
\end{equation*}
Letting $J\to\infty$ in \eqref{eq:volume-sum} therefore gives
$\sum_j e_j<\infty$. On each $I_j\times N$, the trace estimate and \eqref{eq:volume-velocity} give
$\int_{I_j}|z'|\le C_L(H_j+e_j)$. Consequently
\begin{equation}
 \label{eq:volume-summable}
 \sum_{j=0}^\infty Q_j<\infty,
 \qquad \int_a^\infty|z'(r)|\,dr<\infty.
\end{equation}
Thus $z(r)$ converges. The strict bound $|z|<\rho/2$ proved above implies that its limit $z_\infty$ satisfies  $|z_\infty|\le \rho/2$.

Put $r_j=a+jL$. Theorem~\ref{thm:cmc-foliation} bounds $u$ uniformly, and $H(\infty)=0$ together with \eqref{eq:volume-concavity} gives
\begin{equation*}
 H_j=\int_{r_j}^\infty\fint_N u|A|^2\,dV_\gamma\,dr.
\end{equation*}
For $r\in I_k$, \eqref{eq:volume-velocity} and the trace estimate \eqref{eq:trace-estimate_plus} give $\|A(r)\|_{L^2(N)}\le C Q_k$. Uniform equivalence of the leaf measures yields
\begin{equation*}
 H_j\le C\sum_{k=j}^\infty\int_{I_k}\|A(r)\|_{L^2(N)}^2\,dr\le C_L\sum_{k=j}^\infty Q_k^2.
\end{equation*}
For fixed $j$, let $J\to\infty$ in \eqref{eq:volume-transverse}. Since $e_J\to0$, the term $e_J\theta^{J-j}$ tends to zero. Moreover, \eqref{eq:volume-summable} implies convergence of the finite sum on the right hand side of \eqref{eq:volume-transverse} to
$\sum_{k=0}^{\infty}\theta^{|j-k|}Q_k$.
Adding the bound
$H_j\le C_L\sum_{k=j}^{\infty}Q_k^2$ gives
\begin{equation}
 \label{eq:volume-rate-inequality}
 Q_j\le C e_0\theta^j+C\varepsilon\sum_{k=0}^\infty\theta^{|j-k|}Q_k+C\sum_{k=j}^\infty Q_k^2.
\end{equation}

By \eqref{eq:volume-summable}, choose $j_0$ sufficiently large that $\sum_{j\ge j_0}Q_j\le\eta$, and replace $a$ by $a+j_0L$ and reindex the $I_j\times N$ starting at $j_0$. Apply \eqref{eq:volume-transverse} to finite cylinders with this new left endpoint and then let the right endpoint tend to infinity. This gives \eqref{eq:volume-rate-inequality}. The constants are unchanged because the estimate is uniform in the initial parameter. This only changes the eventual multiplicative constant in the decay estimate. Write $M_j=\sup_{k\ge j}Q_k$. For $l\ge j$, split the convolution in \eqref{eq:volume-rate-inequality} into $k<j$, $j\le k<l$, and $k\ge l$. Then
\begin{align*}
 \sum_{k<j}\theta^{l-k}Q_k
   &\le\sum_{k<j}\theta^{j-k}M_k,\\
 \sum_{j\le k<l}\theta^{l-k}Q_k
   &\le\frac{\theta}{1-\theta}M_j,\\
 \sum_{k\ge l}\theta^{k-l}Q_k&\le\frac1{1-\theta}M_j,\quad \sum_{k\ge l}Q_k^2\le\eta M_j.
\end{align*}
Taking the supremum over $l\ge j$ gives
\begin{equation*}
 M_j\le Ce_0\theta^j+C\varepsilon\sum_{k<j}\theta^{j-k}M_k+\left(C\varepsilon\frac{1+\theta}{1-\theta}+C\eta\right)M_j.
\end{equation*}
Choose first $\varepsilon$ and then $\eta$ so that the last coefficient is at most $1/2$. After absorption,
\begin{equation*}
M_j\le2Ce_0\theta^j+2C\varepsilon\sum_{k<j}\theta^{j-k}M_k.
\end{equation*}
Dividing by $\theta^j$ gives
\begin{equation*}
 \frac{M_j}{\theta^j}
 \le2Ce_0+2C\varepsilon\sum_{k<j}\frac{M_k}{\theta^k}.
\end{equation*}
Induction in this last inequality yields $M_j/\theta^j\le2Ce_0(1+2C\varepsilon)^j$, and hence
\begin{equation}
 \label{eq:volume-exponential}
 M_j\le2Ce_0\bigl[\theta(1+2C\varepsilon)\bigr]^j.
\end{equation}
Our earlier choice makes $\theta(1+2C\varepsilon)<1$. This proves exponential decay. Equation~\eqref{eq:volume-velocity} and integration of $z'$ then give exponential decay of $z-z_\infty$ as well. The same equation gives exponential decay of $\|h_-\|_{\mathsf{H}^1}$. The trace estimate \eqref{eq:trace-estimate_plus} gives it for $\|h_+\|_{\mathsf{H}^1}$.
Thus the entire transverse remainder decays exponentially.

The final step is to combine the decay of $\gamma-g_{z_\infty}$ with bounds for $u-1$ and $X$ to establish decay of $G-(dr^2+g_{z_\infty})$.
Lemma~\ref{lem:quantitative-shift} and the uniform bound for $u$ give $\|X(r)\|_{\mathsf{H}^1}\le C\|A(r)\|_{L^2}$.
Equation~\eqref{eq:transverse-normal-speed-bound} gives $\|u-1\|_{\mathsf{H}^2}\le C\|A\|_{L^2}$. Thus $X$ and $u-1$ decay exponentially. Since $\gamma-g_{z_\infty}=h+(g_z-g_{z_\infty})$, the decay of $h$ and $z-z_\infty$ shows that, for $U=G-(dr^2+g_{z_\infty})$, the $L^2$ norm of $U$ on every $[r,r+1]\times N$ decays exponentially.

By Proposition~\ref{prop:slice-path}, in particular \eqref{eq:slice-block-smallness}, the fields $h$, $\partial_r\gamma$, $u-1$, and $X$ have uniformly bounded derivatives of every fixed order on every $I_j\times N$. Derivatives of $\gamma=g_z+h$ along $N$ are uniformly bounded by Lemma~\ref{lem:ricci-flat-decomposition}, since $z(r)$ remains in a compact subset. The metric formula \eqref{eq:slice-metric} gives, for every integer $s\ge0$,
\[
 \sup_j
 \|U\|_{\mathsf{H}^{2s}(\widetilde I_j\times N)}
 \le C_s,
\]
where the enlarged intervals $\widetilde I_j$ have fixed width and $C_s$ is independent of $j$. Combining this uniform bound with the exponential $L^2$ decay established above, interior interpolation on $I_j\Subset\widetilde I_j$ gives
\begin{equation*}
 \|U\|_{\mathsf{H}^s(I_j\times N)}
 \le C_s\|U\|_{L^2(\widetilde I_j\times N)}^{1/2}\|U\|_{\mathsf{H}^{2s}(\widetilde I_j\times N)}^{1/2}
 \le C_s e^{-\beta r_j/2},
\end{equation*}
where $I_j$ and $\widetilde I_j$ have fixed lengths. Given $k$ and $\alpha$, choose $s>k+\alpha+n/2$ and apply Sobolev embedding. After renaming $\beta/2$ as $\beta$, this proves
\begin{equation*}
 G-(dr^2+g_{z_\infty})=O(e^{-\beta r})
\end{equation*}
in $C^{k,\alpha}([r,r+1]\times N)$ for every fixed $k$.

Every translation limit of this end is now $(\mathbb R\times N,dr^2+g_{z_\infty})$.  Along the initially prescribed sequence that limit was $(\mathbb R\times N,dr^2+g_N)$, so the two products are isometric. An isometry between products with compact connected cross sections preserves their slices up to translation and reversal of the $\mathbb R$ coordinate. The pullback of the target radial coordinate has parallel gradient and therefore is $\pm r+c$. Their cross sections are consequently isometric. This finishes the proof of both uniqueness and convergence speed \eqref{eq:speed}.
\end{proof}

\subsection{Optimal decay rate}
\label{sec:optimal-rate}

Once some exponential rate is known, the linearized equation at the limiting metric gives a sharp rate, by a bootstrap argument analogous to the conical argument in \cite{Cheeger-Tian1994}*{Theorem 5.78}. Put $g_\infty=g_{z_\infty}$. On the sufficiently far end, apply Lemma~\ref{lem:fixed-affine-slice} with background $g_\infty$ and normalize the shift by $X\perp\operatorname{Kill}(g_\infty)$. The leafwise coordinate changes $\phi_r$ can be chosen to tend to the identity. Retaining the notation $G,u,\gamma,X$ in these coordinates, set
\[
 q=\gamma-g_\infty,\qquad \delta_{g_\infty}q=0.
\]
All decay statements below include derivatives of every fixed order on $[r,r+1]\times N$.

Suppose $G-(dr^2+g_\infty)=O(e^{-br})$ for some $b>0$. The identity for $H'$ in \eqref{eq:volume-concavity}, with $H(\infty)=0$, and the normal-speed and shift equations used in Lemma~\ref{lem:transverse-preliminary} give
\[
 H=O(e^{-2br}),\qquad u-1=O(e^{-2br}),\qquad X=O(e^{-2br}).
\]
For the shift estimate, the coefficient differences are now $O(e^{-br})$, while $\gamma'=O(e^{-br})$. Taking the $g_\infty$ trace of $\gamma'=2uA+\mathcal L_X\gamma$ and integrating from infinity gives $\operatorname{tr}_{g_\infty}q=O(e^{-2br})$. Let $P_{\mathrm{TT}}$ denote the $L^2(g_\infty)$-orthogonal projection onto the space of transverse-traceless tensors
\[
 \mathrm{TT}_{g_\infty}
 =\{v:\delta_{g_\infty}v=0,\ \operatorname{tr}_{g_\infty}v=0\}.
\]
The decomposition of symmetric tensors into direct sum of transverse traceless part, pure trace part and image of $\delta^*_{g_\infty}$ then gives $q-P_{\mathrm{TT}}q=O(e^{-2br})$. Taylor expansion at $g_\infty$, exactly as in Lemma~\ref{lem:transverse-equation}, yields
\[
 (-\partial_r^2+\Delta^L_{g_\infty})P_{\mathrm{TT}}q =O(e^{-2br}).
\]

When the positive spectrum of $\Delta^L_{g_\infty}$ on $\mathrm{TT}_{g_\infty}$ is nonempty, define
\[
 \lambda_{\mathrm{TT}}
 =\min\bigl(\operatorname{Spec}(\Delta^L_{g_\infty}|_{\mathrm{TT}_{g_\infty}})\cap(0,\infty)\bigr).
\]
On an eigentensor with eigenvalue $\lambda$, the radial equation is $-f''+\lambda f=O(e^{-2br})$. Decay at infinity excludes the constant, linear, and oscillatory homogeneous solutions for $\lambda\le0$, so these components are $O(e^{-2br})$. For $\lambda>0$, the decaying homogeneous solution is $e^{-\sqrt\lambda r}$. The corresponding Green kernels, together with interior elliptic estimates, improve the rate to every
\[
 0<b'<\min\{2b,\sqrt{\lambda_{\mathrm{TT}}}\}.
\]
Iteration gives every rate below $\sqrt{\lambda_{\mathrm{TT}}}$. Once $2b>\sqrt{\lambda_{\mathrm{TT}}}$, the forcing decays strictly faster than $e^{-\sqrt{\lambda_{\mathrm{TT}}}}$, and the same kernels give the endpoint estimate
\[
 G-(dr^2+g_\infty)=O(e^{-\beta_*r}).
\]
More precisely, for some $\delta>0$,
\[
 \gamma-g_\infty=e^{-\sqrt{\lambda_{\mathrm{TT}}}r}v
       +O(e^{-(\sqrt{\lambda_{\mathrm{TT}}}+\delta)r}),\qquad
 v\in\mathrm{TT}_{g_\infty},\quad
 \Delta^L_{g_\infty}v=\lambda_{\mathrm{TT}}v.
\]
Thus $\sqrt{\lambda_{\mathrm{TT}}}$ is sharp for this end if $v\ne0$. If $v=0$, the end converges faster. 

\bibliographystyle{amsalpha}
\bibliography{uniqueness}

@preamble { "\newcommand{\noopsort}[1]{} " }

@article {yaulinearvolumegrowth,
    AUTHOR = {Yau, Shing Tung},
     TITLE = {Some function-theoretic properties of complete {R}iemannian
              manifold and their applications to geometry},
   JOURNAL = {Indiana Univ. Math. J.},
  FJOURNAL = {Indiana University Mathematics Journal},
    VOLUME = {25},
      YEAR = {1976},
    NUMBER = {7},
     PAGES = {659--670},
      ISSN = {0022-2518},
   MRCLASS = {32F05 (31C05 53C55)},
  MRNUMBER = {417452},
MRREVIEWER = {Hung-Hsi Wu},
       DOI = {10.1512/iumj.1976.25.25051},
       URL = {https://doi.org/10.1512/iumj.1976.25.25051},
}

@article {CheegerColding96,
    AUTHOR = {Cheeger, Jeff and Colding, Tobias H.},
     TITLE = {Lower bounds on {R}icci curvature and the almost rigidity of
              warped products},
   JOURNAL = {Ann. of Math. (2)},
  FJOURNAL = {Annals of Mathematics. Second Series},
    VOLUME = {144},
      YEAR = {1996},
    NUMBER = {1},
     PAGES = {189--237},
      ISSN = {0003-486X,1939-8980},
   MRCLASS = {53C21 (53C20 53C23)},
  MRNUMBER = {1405949},
MRREVIEWER = {Joseph\ E.\ Borzellino},
       DOI = {10.2307/2118589},
       URL = {https://doi.org/10.2307/2118589},
}

@article {Cheeger-Colding97I,
    AUTHOR = {Cheeger, Jeff and Colding, Tobias H.},
     TITLE = {On the structure of spaces with {R}icci curvature bounded
              below. {I}},
   JOURNAL = {J. Differential Geom.},
  FJOURNAL = {Journal of Differential Geometry},
    VOLUME = {46},
      YEAR = {1997},
    NUMBER = {3},
     PAGES = {406--480},
      ISSN = {0022-040X},
   MRCLASS = {53C21 (53C20)},
  MRNUMBER = {1484888},
MRREVIEWER = {William P. Minicozzi, II},
       URL = {http://projecteuclid.org/euclid.jdg/1214459974},
}

@article {CN11,
    AUTHOR = {Colding, Tobias H. and Naber, Aaron},
     TITLE = {Characterization of tangent cones of noncollapsed limits with
              lower {R}icci bounds and applications},
   JOURNAL = {Geom. Funct. Anal.},
  FJOURNAL = {Geometric and Functional Analysis},
    VOLUME = {23},
      YEAR = {2013},
    NUMBER = {1},
     PAGES = {134--148},
      ISSN = {1016-443X},
   MRCLASS = {53C20 (53C21)},
  MRNUMBER = {3037899},
MRREVIEWER = {Yu Ding},
       DOI = {10.1007/s00039-012-0202-7},
       URL = {https://doi.org/10.1007/s00039-012-0202-7},
}

@article {CN15_codim4,
    AUTHOR = {Cheeger, Jeff and Naber, Aaron},
     TITLE = {{Regularity of {E}instein manifolds and the codimension 4
              conjecture}},
   JOURNAL = {Ann. of Math. (2)},
  FJOURNAL = {Annals of Mathematics. Second Series},
    VOLUME = {182},
      YEAR = {2015},
    NUMBER = {3},
     PAGES = {1093--1165},
      ISSN = {0003-486X},
   MRCLASS = {53C25 (53C23)},
   MRNUMBER = {3418535},
MRREVIEWER = {Luis Guijarro},
       DOI = {10.4007/annals.2015.182.3.5},
       URL = {https://doi.org/10.4007/annals.2015.182.3.5},
}

@article {SormaniMiniVol,
    AUTHOR = {Sormani, Christina},
     TITLE = {Busemann functions on manifolds with lower bounds on {R}icci
              curvature and minimal volume growth},
   JOURNAL = {J. Differential Geom.},
  FJOURNAL = {Journal of Differential Geometry},
    VOLUME = {48},
      YEAR = {1998},
    NUMBER = {3},
     PAGES = {557--585},
      ISSN = {0022-040X,1945-743X},
   MRCLASS = {53C21 (53C20 53C23)},
  MRNUMBER = {1638053},
MRREVIEWER = {Zhongmin\ Shen},
       URL = {http://projecteuclid.org/euclid.jdg/1214460863},
}

@article {SormaniSublinear,
    AUTHOR = {Sormani, Christina},
     TITLE = {The almost rigidity of manifolds with lower bounds on {R}icci
              curvature and minimal volume growth},
   JOURNAL = {Comm. Anal. Geom.},
  FJOURNAL = {Communications in Analysis and Geometry},
    VOLUME = {8},
      YEAR = {2000},
    NUMBER = {1},
     PAGES = {159--212},
      ISSN = {1019-8385,1944-9992},
   MRCLASS = {53C24 (53C21 53C23)},
  MRNUMBER = {1730892},
MRREVIEWER = {William\ P.\ Minicozzi, II},
       DOI = {10.4310/CAG.2000.v8.n1.a6},
       URL = {https://doi.org/10.4310/CAG.2000.v8.n1.a6},
}

@article {Zhu2025,
    AUTHOR = {Zhu, Xingyu},
     TITLE = {On the geometry at infinity of manifolds with linear volume
              growth and nonnegative {R}icci curvature},
   JOURNAL = {Trans. Amer. Math. Soc.},
  FJOURNAL = {Transactions of the American Mathematical Society},
    VOLUME = {378},
      YEAR = {2025},
    NUMBER = {1},
     PAGES = {503--526},
      ISSN = {0002-9947,1088-6850},
   MRCLASS = {53C21 (53C23)},
  MRNUMBER = {4840313},
       DOI = {10.1090/tran/9261},
       URL = {https://doi.org/10.1090/tran/9261},
}

@article {Cheeger-Tian1994,
    AUTHOR = {Cheeger, Jeff and Tian, Gang},
     TITLE = {On the cone structure at infinity of {R}icci flat manifolds
              with {E}uclidean volume growth and quadratic curvature decay},
   JOURNAL = {Invent. Math.},
  FJOURNAL = {Inventiones Mathematicae},
    VOLUME = {118},
      YEAR = {1994},
    NUMBER = {3},
     PAGES = {493--571},
      ISSN = {0020-9910,1432-1297},
   MRCLASS = {53C21 (53C25 53C55)},
  MRNUMBER = {1296356},
MRREVIEWER = {Xiao\ Wei\ Peng},
       DOI = {10.1007/BF01231543},
       URL = {https://doi.org/10.1007/BF01231543},
}

@article{ColdingMinicozziUniqueness,
    AUTHOR = {Colding, Tobias Holck and Minicozzi, II, William P.},
     TITLE = {On uniqueness of tangent cones for {E}instein manifolds},
   JOURNAL = {Invent. Math.},
  FJOURNAL = {Inventiones Mathematicae},
    VOLUME = {196},
      YEAR = {2014},
    NUMBER = {3},
     PAGES = {515--588},
      ISSN = {0020-9910,1432-1297},
   MRCLASS = {53C25 (53C21 53C23)},
  MRNUMBER = {3211041},
MRREVIEWER = {Megan\ M.\ Kerr},
       DOI = {10.1007/s00222-013-0474-z},
       URL = {https://doi.org/10.1007/s00222-013-0474-z},
}

@article {HaskinsHeinNordstrom,
    AUTHOR = {Haskins, Mark and Hein, Hans-Joachim and Nordstr\"om,
              Johannes},
     TITLE = {Asymptotically cylindrical {C}alabi-{Y}au manifolds},
   JOURNAL = {J. Differential Geom.},
  FJOURNAL = {Journal of Differential Geometry},
    VOLUME = {101},
      YEAR = {2015},
    NUMBER = {2},
     PAGES = {213--265},
      ISSN = {0022-040X,1945-743X},
   MRCLASS = {53C55 (32Q25 53C25)},
  MRNUMBER = {3399097},
MRREVIEWER = {Yalong\ Shi},
       URL = {http://projecteuclid.org/euclid.jdg/1442364651},
}

@article {NordstromG2,
    AUTHOR = {Nordstr\"om, Johannes},
     TITLE = {Deformations of asymptotically cylindrical {$G_2$}-manifolds},
   JOURNAL = {Math. Proc. Cambridge Philos. Soc.},
  FJOURNAL = {Mathematical Proceedings of the Cambridge Philosophical
              Society},
    VOLUME = {145},
      YEAR = {2008},
    NUMBER = {2},
     PAGES = {311--348},
      ISSN = {0305-0041,1469-8064},
   MRCLASS = {53C29},
  MRNUMBER = {2442130},
MRREVIEWER = {Francisco\ Mart\'in Cabrera},
       DOI = {10.1017/S0305004108001333},
       URL = {https://doi.org/10.1017/S0305004108001333},
}

@article {JiangNaber,
    AUTHOR = {Jiang, Wenshuai and Naber, Aaron},
     TITLE = {{$L^2$} curvature bounds on manifolds with bounded {R}icci
              curvature},
   JOURNAL = {Ann. of Math. (2)},
  FJOURNAL = {Annals of Mathematics. Second Series},
    VOLUME = {193},
      YEAR = {2021},
    NUMBER = {1},
     PAGES = {107--222},
      ISSN = {0003-486X,1939-8980},
   MRCLASS = {53B20 (35A21 53C21)},
  MRNUMBER = {4199730},
MRREVIEWER = {Nan\ Li},
       DOI = {10.4007/annals.2021.193.1.2},
       URL = {https://doi.org/10.4007/annals.2021.193.1.2},
}

@article {AndersonConvergence,
    AUTHOR = {Anderson, Michael T.},
     TITLE = {Convergence and rigidity of manifolds under {R}icci curvature
              bounds},
   JOURNAL = {Invent. Math.},
  FJOURNAL = {Inventiones Mathematicae},
    VOLUME = {102},
      YEAR = {1990},
    NUMBER = {2},
     PAGES = {429--445},
      ISSN = {0020-9910,1432-1297},
   MRCLASS = {53C23 (53C21 58D27)},
  MRNUMBER = {1074481},
MRREVIEWER = {Gudlaugur\ Thorbergsson},
       DOI = {10.1007/BF01233434},
       URL = {https://doi.org/10.1007/BF01233434},
}

@article {ChenChen,
    AUTHOR = {Chen, Gao and Chen, Xiuxiong},
     TITLE = {Gravitational instantons with faster than quadratic curvature
              decay. {I}},
   JOURNAL = {Acta Math.},
  FJOURNAL = {Acta Mathematica},
    VOLUME = {227},
      YEAR = {2021},
    NUMBER = {2},
     PAGES = {263--307},
      ISSN = {0001-5962,1871-2509},
   MRCLASS = {53C26 (14D21 32L25 53C80 83C05)},
  MRNUMBER = {4366415},
MRREVIEWER = {Jason\ Dean\ Lotay},
       DOI = {10.4310/acta.2021.v227.n2.a2},
       URL = {https://doi.org/10.4310/acta.2021.v227.n2.a2},
}

@article {DaiWangWei,
    AUTHOR = {Dai, Xianzhe and Wang, Xiaodong and Wei, Guofang},
     TITLE = {On the stability of {R}iemannian manifold with parallel
              spinors},
   JOURNAL = {Invent. Math.},
  FJOURNAL = {Inventiones Mathematicae},
    VOLUME = {161},
      YEAR = {2005},
    NUMBER = {1},
     PAGES = {151--176},
      ISSN = {0020-9910,1432-1297},
   MRCLASS = {53C27 (58J60)},
  MRNUMBER = {2178660},
MRREVIEWER = {Frederik\ Witt},
       DOI = {10.1007/s00222-004-0424-x},
       URL = {https://doi.org/10.1007/s00222-004-0424-x},
}

@book {Besse,
    AUTHOR = {Besse, Arthur L.},
     TITLE = {Einstein manifolds},
    SERIES = {Classics in Mathematics},
      NOTE = {Reprint of the 1987 edition},
 PUBLISHER = {Springer-Verlag, Berlin},
      YEAR = {2008},
     PAGES = {xii+516},
      ISBN = {978-3-540-74120-6},
   MRCLASS = {53C25 (53-02)},
  MRNUMBER = {2371700},
}

@article {SiomonAsymptotics,
    AUTHOR = {Simon, Leon},
     TITLE = {Asymptotics for a class of nonlinear evolution equations, with
              applications to geometric problems},
   JOURNAL = {Ann. of Math. (2)},
  FJOURNAL = {Annals of Mathematics. Second Series},
    VOLUME = {118},
      YEAR = {1983},
    NUMBER = {3},
     PAGES = {525--571},
      ISSN = {0003-486X,1939-8980},
   MRCLASS = {58G11 (35B40 49F99 58E20)},
  MRNUMBER = {727703},
MRREVIEWER = {Helmut\ Kaul},
       DOI = {10.2307/2006981},
       URL = {https://doi.org/10.2307/2006981},
}

@article {AndersonEinsteinMetric,
    AUTHOR = {Anderson, Michael T.},
     TITLE = {Ricci curvature bounds and {E}instein metrics on compact
              manifolds},
   JOURNAL = {J. Amer. Math. Soc.},
  FJOURNAL = {Journal of the American Mathematical Society},
    VOLUME = {2},
      YEAR = {1989},
    NUMBER = {3},
     PAGES = {455--490},
      ISSN = {0894-0347,1088-6834},
   MRCLASS = {53C20 (53C25 58D17 58G30)},
  MRNUMBER = {999661},
MRREVIEWER = {Maung\ Min-Oo},
       DOI = {10.2307/1990939},
       URL = {https://doi.org/10.2307/1990939},
}

@phdthesis{haslhoferPhD,
  author = {Haslhofer, Robert},
  title = {Stability of Ricci-flat Spaces and
Singularities in 4d Ricci Flow},
  school = {ETH Zurich},
  year = {2012},
  type = {PhD dissertation},
   URL = {https://www.math.toronto.edu/roberth/papers/Haslhofer_thesis.pdf},
}

@article {SunZhang,
    AUTHOR = {Sun, Song and Zhang, Junsheng},
     TITLE = {No semistability at infinity for {C}alabi-{Y}au metrics
              asymptotic to cones},
   JOURNAL = {Invent. Math.},
  FJOURNAL = {Inventiones Mathematicae},
    VOLUME = {233},
      YEAR = {2023},
    NUMBER = {1},
     PAGES = {461--494},
      ISSN = {0020-9910,1432-1297},
   MRCLASS = {53C55 (32Q26)},
  MRNUMBER = {4602001},
MRREVIEWER = {Antonella\ Nannicini},
       DOI = {10.1007/s00222-023-01187-4},
       URL = {https://doi.org/10.1007/s00222-023-01187-4},
}

@misc{BamlerLai2025,
  author        = {Bamler, Richard H. and Lai, Yi},
  title         = {The {PDE-ODI} principle and cylindrical mean curvature flows},
  year          = {2025},
  eprint        = {2512.25050},
  archivePrefix = {arXiv},
  primaryClass  = {math.DG},
  doi           = {10.48550/arXiv.2512.25050},
  url           = {https://arxiv.org/abs/2512.25050}
}

\end{document}